\documentclass[12pt,titlepage]{article}

\usepackage{graphicx,color}
\usepackage{amsmath,amssymb,amsthm}
\renewcommand{\qed}{\hfill\blacksquare}
\newtheorem{thm}{Theorem}
\newtheorem{lem}[thm]{Lemma}
\newtheorem{cor}[thm]{Corollary}

\newtheorem{prop}[thm]{Proposition}

\newtheorem{defi}{Definition}

\newtheorem{remark}{Remark}

\def\beq{\begin{equation}}\def\eeq{\end{equation}}
\def\beqn{\begin{eqnarray}}\def\eeqn{\end{eqnarray}}

\def\qed{\ifhmode\unskip\nobreak\fi\quad\ifmmode\Box\else$\Box$\fi}

\newcommand{\thet}{\overline{\vartheta}}

\newcommand\mbf[1]{\mbox{\boldmath$#1$}}
\newcommand\msbf[1]{\mbox{\boldmath\scriptsize$#1$}}

\title{Shannon capacity, Lov\'asz theta number and the Mycielski construction\footnote{Some of the results in this paper were presented at the Szeged Workshop on Discrete Structures (SWORDS2023) held in Szeged, Hungary, in November 2023.}}
\author{\hfil Bence Csonka\thanks
  {Department of Computer Science and Information Theory,
Faculty of Electrical Engineering and Informatics,
Budapest University of Technology and Economics and
MTA-BME Lend\"ulet Arithmetic Combinatorics Research Group, MKH, Budapest, Hungary. Research partially supported by the “TKP2020,
National Challenges Program” of the National Research Development and
Innovation Office (BME NC TKP2020) and also by the Ministry of Innovation and
Technology and the National Research, Development and Innovation
Office within the Artificial Intelligence National Laboratory of Hungary; email: {\tt slenhortag@gmail.com}
 }
\and G\'abor Simonyi\thanks{Alfr\'ed R\'enyi Institute of Mathematics (HUN-REN), Budapest, Hungary and
Department of Computer Science and Information Theory,
Faculty of Electrical Engineering and Informatics,
Budapest University of Technology and Economics. Research partially supported by the National Research, Development and Innovation Office (NKFIH) grants K--132696 and SNN-135643 of NKFIH Hungary; email: {\tt simonyi@renyi.hu}
}
}

\date{}

\begin{document}

\maketitle

\begin{abstract}
We investigate the effect of the well-known Mycielski construction on the Shannon capacity of graphs and on one of its most prominent upper bounds, the (complementary) Lov\'asz theta number.
We prove that if the Shannon capacity of a graph, the distinguishability graph of a noisy channel, is attained by some finite power, then its Mycielskian has strictly larger Shannon capacity than the graph itself. For the complementary Lov\'asz theta function we show that its value on the Mycielskian of a graph is completely determined by its value on the original graph, a phenomenon similar to the one discovered for the fractional chromatic number by Larsen, Propp and Ullman. We also consider the possibility of generalizing our results on the Sperner capacity of directed graphs and on the generalized Mycielsky construction. Possible connections with what Zuiddam calls the asymptotic spectrum of graphs are discussed as well.



\end{abstract}

\section{Introduction}

\bigskip
\par\noindent
The zero-error capacity of a noisy channel was introduced by Shannon \cite{Sha56} who showed that it can be expressed as an asymptotic parameter of a graph $G$.
If $G$ is the distinguishability graph of a noisy channel, that is, its vertex set is the input alphabet of the channel and two letters form an edge if they cannot result in the same output letter, then the maximum number of pairwise ditinguishable input letters is equal to $\omega(G)$, the clique number of graph $G$. The distinguishability relation is extended to sequences of letters by the ${\rm OR}$-product of graphs which is defined as follows.

\medskip
\par\noindent
\begin{defi} \label{defi:ORprod}
For two graphs $F$ and $G$ their OR-product $F\cdot G$ is defined by
$$V(F\cdot G)=V(F)\times V(G)$$ and
$$E(F\cdot G)=\{\{(f,g),(f',g')\}: f,f'\in V(F), g,g'\in V(G),$$$$\{f,f'\}\in E(F)\ {\rm or}\ \{g,g'\}\in E(G)\}.$$
The $t^{\rm th}$ OR-power $G^t$ of a graph $G$ is the $t$-fold OR-product of $G$ with itself.
\end{defi}

\medskip
\par\noindent
Since two sequences of length $t$ consisting of the input letters cannot result in the same output sequence exactly if in at least one coordinate they cannot result in the same output letter, the maximum number of pairwise distinguishable sequences of length $t$ is given by $\omega(G^t)$. Thus the (logarithmic) zero-error capacity of the channel with distinguishability graph $G$ is $c_{\rm OR}(G):=\lim_{t\to\infty}\frac{1}{t}\log_2\omega(G^t)$. It is somewhat more convenient for our discussion if we do not have to normalize all our bounds by taking the logarithm, so we will use the non-logarithmic version.

\medskip
\par\noindent
\begin{defi}\label{defi:orcap}
  The (non-logarithmic) Shannon OR-capacity of graph $G$ is defined by
  $$C_{\rm OR}(G):=\lim_{t\to\infty}\sqrt[t]{\omega(G^t)}.$$
\end{defi}

\medskip
\par\noindent
It is easy to see that for every graph $G$ and any fixed $k$ we have $\sqrt[k]{\omega(G^k)}\le C_{\rm OR}(G)$. This relation already implies the well-known fact that the above limit always exists by Fekete's Lemma.
We remark that Shannon~\cite{Sha56} used a complementary language when defining graph capacity, modelling the channel by its confusability graph which is just the complement of the distinguishability graph. The above definition is, however, more adequate for our discussion, just as it is the case in Chapter 11 of the book \cite{CsK}. Nevertheless, we use the ``OR'' subscript in the above definition to signify this difference.

\medskip
\par\noindent
It is easy to show that Shannon OR-capacity satisfies
$$\omega(G)\le C_{\rm OR}(G)\le \chi(G),$$ where $\chi(G)$ denotes the chromatic number of the graph.

Thus if $\chi(G)=\omega(G)$, then $C_{\rm OR}(G)$ also shares their value. The smallest graph for which $\chi(G)>\omega(G)$ is $C_5$, the cycle of length $5$. Shannon~\cite{Sha56} has shown that its capacity value is at least $\sqrt{5}$ by presenting a clique of size $5$ in its second OR-power (though in a complementary language, as remarked above). The matching upper bound was proven only more than two decades later by Lov\'asz~\cite{LL79} as a first application of his $\vartheta$-number introduced in \cite{LL79}. The value of Shannon OR-capacity is still unkonwn for longer odd cycles and their complements.

\medskip
\par\noindent
It is well-known that the gap between $\chi(G)$ and $\omega(G)$ can be arbitrarily large and one of the best-known constructions producing a sequence of graphs with that gap going to infinity is due to Mycielski \cite{Myc}. This construction works iteratively. For any graph $G$ it produces a graph $M(G)$, called its Mycielskian, with the property that $\omega(M(G))=\omega(G)$ while
$\chi(M(G))=\chi(G)+1$.

\medskip
\par\noindent
\begin{defi}\label{defi:Myc}
  Let $G$ be a simple graph.
  Its Mycielskian $M(G)$ is defined on the vertex set $$V(M(G))=V(G)\times\{0,1\}\cup \{z_{M(G)}\}$$
  with edge set $$E(M(G))=\{\{(v,0)(w,i)\}: \{v,w\}\in E(G), i\in\{0,1\}\}\cup\{\{z_{M(G)},(v,1)\}: v\in V(G)\}.$$

\end{defi}

\smallskip
\par\noindent
Note that $C_5$ is isomorphic to $M(K_2)$, where $K_n$ is the complete graph on $n$ vertices.

\medskip
\par\noindent
In this paper we investigate the effect of the Mycielski construction on the value of Shannon OR-capacity and on one of its most prominent upper bounds, the complementary Lov\'asz theta number.
For the former we will prove the following theorem.

\medskip
\par\noindent
\begin{thm}\label{thm:fintcap}
  If $G$ is a graph that attains its Shannon capacity in finite length, i.e., there exists some positive integer $k$ for which $C_{\rm OR}(G)=\sqrt[k]{\omega(G^k)}$, then
  $$C_{\rm OR}(M(G))>C_{\rm OR}(G).$$
\end{thm}

\medskip
\par\noindent
Larsen, Propp and Ullman \cite{LPU} proved the notable fact that the fractional chromatic number $\chi_f(M(G))$ of the Mycielskian of a graph $G$ is determined by the value of the fractional chromatic number of the original graph. In fact, they showed the validity of the neat formula $$\chi_f(M(G))=\chi_f(G)+\frac{1}{\chi_f(G)}.$$
For the definition and basic properties of the fractional chromatic number we refer to \cite{SchU}. Here we only note that $\chi_f(G)$ is also an upper bound of $C_{\rm OR}(G)$ that satisfies $\chi_f(G)\le\chi(G)$, similarly to the complementary Lov\'asz theta number $\thet(G)$ (which actually also satisfies $\thet(G)\le\chi_f(G)$ for every graph $G$ as proven by Lov\'asz \cite{LL79}) that we will introduce in the next section. Our main result on the latter will be that
$\thet(M(G))$ is also determined by the value of $\thet(G)$. Although more complicated than in case of the fractional chromatic number, we also provide a formula giving this dependence. In particular, we will prove the following statement.

\begin{thm}\label{main}
For every nonempty graph $G$
\[
\thet\left(M(G)\right) = \frac43\thet(G)\cos{\left(\frac13{\rm arccos}\left(1 - \frac{27}{4\thet(G)}+\frac{27}{4\thet^2(G)}\right)\right)}-\frac13\thet(G)+1.
\]
\end{thm}

\begin{remark}
  It is worth noting the case $\thet(G) = 3$ which gives us a simple expression, for example, for $G=K_3$ we get \newline $$\thet(M(K_3)) = 4\cos{\left(\frac{2\pi}{9}\right)}.$$
\end{remark}

\medskip
\par\noindent
The paper is organized as follows. In the next two sections we give the proof of Theorem~\ref{thm:fintcap} and Theorem~\ref{main}, respectively. In Section~\ref{sect:egyebek} we discuss some further related problems, in particular about Sperner capacity and about the so-called generalized Mycielski construction. We also elaborate on the question whether our Theorem~\ref{main} may have anything to do with the complementary Lov\'asz theta number belonging to what Zuiddam~\cite{Zuiddam} calls the asymptotic spectrum of graphs. The paper is concluded with some open problems.

\section{Shannon capacity of the Mycielskian}

\medskip
\par\noindent
The main observation needed for proving Theorem~\ref{thm:fintcap} is the following Lemma.

\medskip
\par\noindent
\begin{lem}~\label{lem:felemel}
There is a clique of size $n^n$ in $[M(K_n)\setminus\{z_{M(K_n)}\}]^n$ every vertex of which has a coordinate belonging to the set $V(K_n)\times\{1\}$.
\end{lem}

\medskip
\par\noindent
\proof
Let the vertices of $K_n$ be $0,1,\dots,n-1$.
Consider the set $B_n:=\{(0,0),(1,0),\dots,(n-1,0)\}^n$, that is the set of all sequences of length $n$ formed by vertices of the $V(K_n)\times \{0\}$ part in our graph $M(K_n)$. Clearly, these sequences form a complete subgraph in $[M(K_n)]^n$.

\medskip
\par\noindent
Define
$B_n^{(j)}\subseteq B_n$ for all $j\in\{0,\dots,n-1\}$ as follows:
$$B_n^{(j)}:=\{(x_1,0)(x_2,0)\dots (x_n,0)\in B_n: \sum_{i=1}^n x_i\equiv j\pmod n\}.$$
Then let $\hat{B}_n^{(0)}= \{(x_1,1)(x_2,0)\dots (x_n,0): (x_1,0)\dots (x_n,0)\in B_n^{(0)}\}$ and in general
$$\hat{B}_n^{(j)}:=\{(x_1,0)\dots (x_j,0)(x_{j+1},1)(x_{j+2},0)\dots (x_n,0): (x_1,0)\dots (x_n,0)\in B_n^{(j)}\}.$$
Furthermore, we let
$$\hat{B}_n:=\hat{B}_n^{(0)}\cup\dots\cup\hat{B}_n^{(n-1)}.$$
We claim that $\hat{B}_n$ induces a clique of size $n^n$ in $[M(K_n)]^n$. Indeed,
if $(x_1,0)\dots (x_{h+1},1)\dots (x_n,0)\in \hat{B}_n^{(h)}, (y_1,0)\dots (y_{r+1,}1)\dots (y_n,0)\in \hat{B}_n^{(r)}$ for $h\neq r$ (that is, the coordinates where we ``lifted'' the corresponding entry of two sequences from $V(K_n)\times\{0\}$ to $V(K_n)\times\{1\}$ are different), then the two new sequences are adjacent in the same coordinate where they were adjacent before the ``lifting''. If, on the other hand, we consider two sequences from the same subset $\hat{B}_n^{(j)}$ like
$(x_1,0)\dots (x_{j+1},1)\dots (x_n,0)$ and $(y_1,0)\dots (y_{j+1},1)\dots (y_n,0)$, then there must be some $i\neq j+1$ where $x_i\neq y_i$ and therefore these two sequences are also adjacent in $[M(K_n)]^n$. The latter is true because the sequences $(x_1,0)\dots (x_{j+1},0)\dots (x_n,0), (y_1,0)\dots (y_{j+1},0)\dots (y_n,0)\in B_n^{(j)}$ must have at least two coordinates where they differ, otherwise we could not have
$\sum_{i=1}^n x_i\equiv \sum_{i=1}^n y_i\pmod n$. So even if one of these two coordinates was the $(j+1)^{\rm th}$ (and thus $(x_{j+1},1)$ and $(y_{j+1},1)$ are not adjacent any more), the two sequences will still be adjacent because of the other coordinate where they differed.

\medskip
\par\noindent
So $\hat{B}_n$ induces a clique of size $n^n$ in $[M(K_n)]^n$ and the sequences forming the vertices in it all have a ``lifted'' coordinate, that is one representing a vertex in the $V((K_n)\times\{1\}$ part of the vertex set of $M(K_n)$.
This completes the proof of the Lemma.
\qed

\medskip
\par\noindent
Note that the above lemma has the following immediate consequence, a special case of Theorem~\ref{thm:fintcap}, which may be of some interest on its own right.

\medskip
\par\noindent
\begin{cor}\label{lem:Kn}
  For the complete graph $K_n$ on $n$ vertices we have
  $$C_{\rm OR}(M(K_n))\ge \sqrt[n]{n^n+1}>n.$$
 \end{cor}

\proof
Consider the clique of size $n^n$ in $[M(K_n)\setminus\{z_{M(K_n)}\}]^n$ whose existence is shown in Lemma~\ref{lem:felemel}. Since all the sequences that are vertices of this clique contain a coordinate from $V(K_n)\times \{1\}=N(z_{M(K_n)})$ (where $N(v)$ stands for the neighborhood of vertex $v$), the all-$z_{M(K_n)}$ sequence of length $n$ is adjacent to all vertices of this clique in $[M(K_n)]^n$. So adding this sequence to our clique of size $n^n$ we obtain a clique with $n^n+1$ vertices, thus showing $\omega([M(K_n)]^n)\ge n^n+1$. This implies $C_{\rm OR}(M(K_n))\ge\sqrt[n]{n^n+1}$ as stated.
\qed

\medskip
\par\noindent
We remark that for $n=2$ this gives $C_{\rm OR}(C_5)=C_{\rm OR}(M(K_2))\ge\sqrt{5}$, which is exactly Shannon's lower bound for the capacity of the $5$-cycle that Lov\'asz \cite{LL79} proved to be tight.

\medskip
\par\noindent
To complete the proof of Theorem~\ref{thm:fintcap} the following simple lemma will be useful.
\medskip
\par\noindent
\begin{lem} \label{lem:benne}
For every graph $G$ and positive integer $t$ we have
$$M(G^t)\subseteq [M(G)]^t.$$
\end{lem}

\proof
Using the notation ${\mbf v}=v_1v_2\dots v_t\in V(G^t)$
the following function $\varphi$ is an embedding of $M(G^t)$ into $[M(G)]^t$:

$$\varphi: ({\mbf v},h)\mapsto (v_1,h)(v_2,h)\dots (v_t,h)\ {\rm for\ } h\in\{0,1\},$$
while $z_{M(G^t)}$ is mapped to the sequence of length $t$ given by $z_{M(G)}\dots z_{M(G)}\in [V(M(G))]^t$.
\medskip
\par\noindent
It is straightforward to check that the image $\varphi(V(M(G^t)))$ of $V(M(G^t))$ induces a subgraph isomorphic to $M(G^t)$ in $[M(G)]^t$.
\qed

\medskip
\par\noindent
{\em Proof of Theorem~\ref{thm:fintcap}.}
Let $G$ be a graph with $C_{\rm OR}(G)=\sqrt[k]{\omega(G^k)}=\sqrt[k]{N}$ where $k\ge 1$ is a positive integer and $N=\omega(G^k)$.
Since $K_N\subseteq G^k$, we also have $M(K_N)\subseteq M(G^k)\subseteq [M(G)]^k,$ where the last relation is by Lemma~\ref{lem:benne}.
This further implies $$[M(K_N)]^N\subseteq [M(G^k)]^N\subseteq [M(G)]^{kN},$$
so by the proof of Corollary~\ref{lem:Kn} we can write
$$N^N+1\le \omega([M(K_N)]^N)\le\omega([M(G^k)]^N)\le\omega([M(G)]^{kN}).$$
The latter implies
$$C_{\rm OR}(M(G))\ge \sqrt[kN]{\omega([M(G)]^{kN})}\ge\sqrt[kN]{N^N+1}>\sqrt[k]{N}=C_{\rm OR}(G)$$
completing the proof.
\qed

\medskip
\par\noindent
\begin{remark}
We note that although the above lower bound is tight for $M(K_2)=C_5$, we have no reason to believe that it would be tight in general. In fact, for $C_{\rm OR}(M(C_5))$ it gives the lower bound $\sqrt[10]{5^5+1}\approx 2.23614$, and for this specific graph the somewhat larger lower bound $\sqrt[4]{28}\approx 2.30036$ is proven in \cite{coletrips}.
\end{remark}

\section{Complementary Lov\'asz $\vartheta$-number of the Mycielskian} \label{sect:theta}

\medskip
\par\noindent
\begin{defi}
    An orthonormal representation of a simple graph $G$ on vertex set $V=\{1,\dots,n\}$ assigns to each $i\in V$ a unit vector $u_i \in \mathbb{R}^d$ (for some appropriate positive integer $d$) such that $\langle u_i,u_j\rangle=0$, whenever $ij \notin E(G)$. An orthonormal representation of the complementary graph is called a dual orthonormal representation.
\end{defi}
\smallskip
\begin{defi}
    The complementary Lov\'asz $\vartheta$-number of a graph $G$ is
    \[
    \thet(G):= \min{\max_{i \in V}\frac{1}{\langle c,u_i\rangle^2}},
    \]
    where the minimum is taken over all dual orthonormal representations $\{u_i : i \in V \}\subseteq \mathbb{R}^d$
of $G$ and all unit vectors $c \in \mathbb{R}^d$. (The dimension $d$ is also chosen so that the minimal value could be attained.)
\end{defi}

In the proof of Theorem \ref{main} we will use some equivalent definitons of $\thet$. First we recall the definition of the invariant called strict vector chromatic number that was defined and shown to be equivalent to the complementary Lov\'asz $\vartheta$-number in \cite{KMS}\footnote{The authors of \cite{KMS} mention in their paper that the equality of the strict vector chromatic number and $\thet(G)$ was proven ``with the help of \'Eva Tardos and David Williamson''.}, for the proof of equivalence see also \cite{LLGG}.

\smallskip
\begin{defi}{\rm (\cite{KMS})}
    A strict vector $t$-coloring of graph $G$ in $\mathbb{R}^d$ assigns to each $i \in V(G)$ a unit vector $u_i \in \mathbb{R}^d$ such that $\langle u_i,u_j\rangle = -\frac{1}{t-1}$, whenever $ij \in E(G)$. The strict vector chromatic number of a graph $G$ is
    \[
    \chi_{v}(G) := \min{\{t\in \mathbb{R}:G\ \text{{\rm admits a strict vector $t$-coloring for some positive integer} $d$}\}}.
    \]
\end{defi}
\smallskip
\begin{thm} {\rm (\cite{KMS})}

    \[
    \thet(G) = \chi_{v}(G)
    \]
    for every graph $G$.
\end{thm}

\medskip
\par\noindent
To give another known expression of $\thet(G)$ we need the following notions.

\smallskip
\begin{defi}
    Let $A \in \mathbb{R}^{n\times n}$ and let ${\rm Sp}(A)$ denote the set of eigenvalues of the matrix $A$. Furthermore, let $\lambda_{{\rm max}}(A)$ and $\lambda_{{\rm min}}(A)$ be the maximal and minimal elements in ${\rm Sp}(A)$, respectively.
\end{defi}

\medskip
\par\noindent
The following formula for $\thet(G)$ is given as Theorem 6 in \cite{LL79}, see also in Proposition 11.9 of \cite{LLGG}.
\smallskip
\begin{thm}\label{spectra} {\rm (Lov\'asz \cite{LL79})}
    For every graph $G$
    \[
    \thet(G) = 1 + \max_T{\frac{\lambda_{{\rm max}}(T)}{\left|\lambda_{{\rm min}}(T)\right|}},
    \]
    where $T \in \mathbb{R}^{n\times n}$ ranges over all symmetric nonzero matrices with $T_{ij} = 0$ for $ij \notin E(G)$.
\end{thm}

\smallskip
\par\noindent
Note that since we deal with simple graphs only, we will always have $T_{ii}=0$ for the matrices feasible for Theorem~\ref{spectra}. This implies that the trace, that is the sum of eigenvalues of these matrices is equal to $0$. This implies that the smallest eigenvalue $\lambda_n$ of the matrix attaining the maximum above is always negative.

\medskip
\par\noindent
We will need the following well-known theorem from linear algebra (cf. e.g. Section 6.2 in \cite{linalg}).

\smallskip
\begin{thm}\label{symmetric}
  Let $A$ be a real symmetric matrix. Then $A$ is orthogonally diagonalizable, that is, $A = P\Lambda P^\top$, where $\Lambda = {\rm diag}(\lambda_1,...,\lambda_n)$ is a real diagonal matrix having the eigenvalues of $A$ in its main diagonal and the $i$th coulumn of $P$ is a unit length eigenvector corresponding to the eigenvalue $\lambda_i$. In this case the eigenvectors give an orthonormal basis of $\mathbb{R}^n$.
\end{thm}

\medskip
\par\noindent
In the proof of Theorem \ref{main} we will have to find the solutions of a cubic equation, therefore the following well-known formula will be useful for us (see e.g. in \cite{wikicube}).

\begin{thm}\label{equation}
    Given the cubic equation $ax^3 + bx^2 + cx + d = 0$ let
    \[
    p: = \frac{3ac - b^2}{3a^2},
    \]
    \[
    q:= \frac{2b^3 - 9abc + 27a^2d}{27a^3}.
    \]
    Then the solutions of the equation can be written as
    \[
    x_k = 2\sqrt{-\frac{p}{3}}\cos{\left(\frac13{\rm arccos}\left(\frac{3q}{2p}\sqrt{-\frac{3}{p}}\right)-\frac{2\pi k}{3}\right)} -\frac{b}{3a} \hspace{0.5 cm}\text{where }k=0,1,2.
    \]
\end{thm}
\bigskip
\par\noindent
    {\em Proof of Theorem \ref{main}.} First we show the upper bound, i.e. the inequality
    \begin{equation}\label{eq:upb}
\thet\left(M(G)\right) \leq \frac43\thet(G)\cos{\left(\frac13{\rm arccos}\left(1 - \frac{27}{4\thet(G)}+\frac{27}{4\thet^2(G)}\right)\right)}-\frac13\thet(G)+1.
\end{equation}
    Let $V(G) = \{1,...,n\}$ , $t:=\chi_v(G)$ and let $\{v_i:i\in V\} \subseteq \mathbb{R}^d$ be an optimal strict vector $t$-coloring of $G$, where $v_i = (v_i^{(1)},..,v_i^{(d)})$. Now we construct a strict vector coloring of $M(G)$ in $\mathbb{R}^{d+1}$. Let the axis of the last coordinate be the $z$-axis.

\medskip
\par
\begin{figure} [h!]
    \centering
    \includegraphics[width = 13 cm]{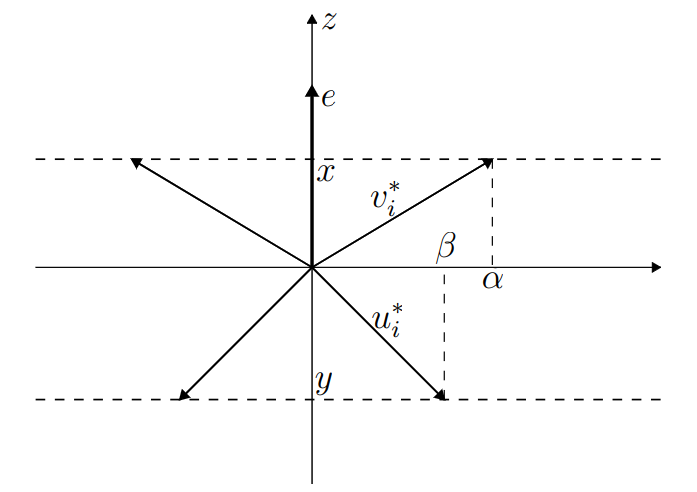}
    \caption{The location of the vectors of our strict vector $\hat{t}$-coloring of $M(G)$ in $\mathbb{R}^{d+1}$. It may be instructive to think about the special case of $G=K_2$ with the optimal strict vector coloring consisting of the two vectors being $1$-dimensional with their unique coordinate being equal to $-1$ and $1$, respectively. Then the five vectors on the picture are in $\mathbb{R}^2$ giving an optimal strict vector coloring of $M(K_2)\cong C_5$ if the values of $\alpha, \beta, x, y$ are chosen appropriately. In the more general case when $G$ has several vertices, we can think about this picture so, that the two dotted lines represent $d$-dimensional hyperplanes. The one defined by  $z=x$ contains all the vectors assigned to vertices $(i,0)$ and the one defined by $z=-y$ contains all the vectors assigned to the ``twin'' vertices $(i,1)$. (Note that $y$ is considered to be a non-negative value, so in the picture it means the distance from the origin.)}
    \label{coordinate}
\end{figure}

\noindent
In our strict vector coloring of $M(G)$ let us assign to each $(i,0)\in V(M(G))$ the vector $v_i^*:=(\alpha v_i^{(1)},...,\alpha v_i^{(d)},x) \in \mathbb{R}^{d+1}$ where the appropriate non-negative real values of $\alpha$ and $x$ are to be chosen later. Similarly, we assign to each ``twin'' vertex $(i,1)$ the vector $u_i^* = (\beta v_i^{(1)},...,\beta v_i^{(d)},-y) \in \mathbb{R}^{d+1}$, where the appropriate non-negative real values of $\beta$ and $y$ are also to be chosen later. Finally, the vector $e := (0,0,...,0,1) \in \mathbb{R}^{d+1}$ is assigned to $z_{M(G)}\in V(M(G))$. (See Figure \ref{coordinate} for an illustration.)

\medskip
\par\noindent
We have to verify that $x,y,\alpha,\beta$ can be chosen so that the system of the vectors $\{v_1^*,...,v_n^*,u_1^*,..,u_n^*,e\}$ gives a strict vector coloring for $M(G)$. This means that these vectors should satisfy the following conditions for some parameter $\hat{t}$:
\begin{enumerate}
    \item For every $i$, $|v_i^*|=|u_i^*|=1$.
    \item If $ij \in E(G)$, then $\langle v_i^*,v_j^*\rangle = -\frac{1}{\hat{t}-1}$.
    \item If $ij \in E(G)$, then $\langle v_i^*,u_j^*\rangle = -\frac{1}{\hat{t}-1}$.
    \item For all $i$, $\langle u_i^*,e \rangle = -\frac{1}{\hat{t}-1}$.
\end{enumerate}
\smallskip
\par\noindent
If these conditions are satisfied then
the given system of vectors forms a strict vector $\hat{t}$-coloring.

\smallskip
\par\noindent
The two equations in 1. give us
\begin{gather*}
1 = |v_i^*|^2 = \sum_{j=1}^{d}\left(\alpha v_i^{(j)}\right)^2 + x^2 = \alpha^2\sum_{j=1}^d\left(v_i^{(j)}\right)^2 + x^2 = \alpha^2 + x^2,\\
1 = |u_i^*|^2 = \sum_{j=1}^{d}\left(\beta v_i^{(j)}\right)^2 + y^2 = \beta^2\sum_{j=1}^d\left(v_i^{(j)}\right)^2 + y^2 = \beta^2 + y^2.
\end{gather*}
The second condition and the fact that the vectors $v_i$ gave a strict vector $t$-coloring of $G$ give that for $ij \in E(G)$ we have
\begin{gather*}
    -\frac{1}{\hat{t}-1} = \langle v_i^*,v_j^*\rangle = \sum_{k=1}^{d}\alpha^2v_i^{(k)}v_j^{(k)} + x^2 =\alpha^2\langle v_i,v_j\rangle + x^2 = -\frac{1}{t-1}\alpha^2 + x^2.
\end{gather*}
Furthermore, for $ij \in E(G)$ the third condition gives
\[
    -\frac{1}{\hat{t}-1} = \langle v_i^*,u_j^*\rangle = \sum_{k=1}^{d}\alpha\beta v_i^{(k)}v_j^{(k)} - xy =-\frac{1}{t-1}\alpha\beta - xy.
\]
Finally, it is straightforward from the fourth condition that $-\frac{1}{\hat{t}-1} = \langle u_i^*,e \rangle = -y$.
To summarize, we have obtained the following system of equations:
\begin{equation}
\centering
\left\{\begin{split}
\alpha^2 + x^2 &= 1 \\
\beta^2 + y^2 &= 1 \\
 \frac{1}{t-1}\alpha^2 - x^2 &=\frac{1}{\hat{t}-1}\\
\frac{1}{t-1}\alpha\beta + xy &=\frac{1}{\hat{t}-1}\\
y  &=\frac{1}{\hat{t}-1}.
\end{split}\right.
\end{equation}

\medskip
\par\noindent
This system of equations provides a solution for $\hat{t}$ that depends only on $t$. After calculation that is presented in Appendix 1, we arrive to the following cubic equation.
\begin{equation}\label{eq:cubic}
\hat{t}^3 + (t-3)\hat{t}^2 + (3-2t-t^2)\hat{t} + (-t^3 +5t^2 -3t-1) = 0.
\end{equation}

\medskip
\par\noindent
According to the formula presented in Theorem \ref{equation} this equation has three solutions that can be given as
\begin{equation}\label{eq:star}
  \hat{t}_k = \frac43t\cos{\left(\frac13{\rm arccos}\left(1 - \frac{27}{4t}+\frac{27}{4t^2}\right)-\frac{2\pi k}{3}\right)}-\frac13t+1 \hspace{0.5 cm} \text{for }k=0,1,2.
\end{equation}
Since $\thet$ is monotone in the sense that $F\subseteq F'\Rightarrow \thet(F)\le \thet(F')$, $G \subseteq M(G)$ implies $\thet(G) \leq \thet(M(G))$. Hence if for some solution $\hat{t}$ of the above equation we have $\hat{t} < \thet(G)$ then $\hat{t}$ cannot be an upper bound on $\thet(M(G))$ and thus it cannot be the value of a strict vector coloring of $M(G)$. The case of $\thet(G) = 1$ is trivial since $\omega(G) \leq \thet(G)$, therefore $G$ is an empty graph or a graph without edges. If $G$ does not have any edge and it is not an empty graph, then $M(G)$ is a star graph with isolated vertices so $\vartheta(M(G)) = 2$. Hence we can assume that $2 \leq \thet(G)$.

\bigskip

\bigskip

\medskip
\par\noindent
We need to show that among the above three solutions of our cubic equation only $\hat{t}_0$ is relevant for us. We can actually show that for $k=1,2$ we have $\hat{t}_k\le 1$ and the above observations already imply that this fact makes those solutions uninteresting for us. Details of this calculation are given in Appendix 2. So we are interested only in $x_0$ that is what we obtain by substituting $k=0$.

\bigskip
\par\noindent
Substituting $k=0$ we obtain the expression on the right hand side of (\ref{eq:upb}), thus the claimed upper bound is proven.

\bigskip
\par\noindent
Now we have to show the reverse inequality, that is
\begin{equation}\label{eq:lowb}
\frac43\thet(G)\cos{\left(\frac13{\rm arccos}\left(1 - \frac{27}{4\thet(G)}+\frac{27}{4\thet^2(G)}\right)\right)}-\frac13\thet(G)+1 \leq \thet\left(M(G)\right).
\end{equation}
The left hand side here in (\ref{eq:lowb}) is a function of $t=\thet(G)$ that we will denote by $m(t)$.
To prove (\ref{eq:lowb}) we will use Theorem \ref{spectra}. Let $t:=\thet(G)$ and assume that $T \in \mathbb{R}^{n \times n}$ is a matrix obtaining the maximum in said theorem for graph $G$, that is, $T$ is a symmetric matrix with $T_{ij} = 0$ for $ij \notin E(G)$ such that $t = 1 + \frac{\lambda_{{\rm max}}(T)}{|\lambda_{{\min}}(T)|}$. Let $\lambda_1\geq \lambda_2 \geq...\geq \lambda_n$ be the eigenvalues of $T$ and let the corresponding orthonormal system of eigenvectors be $\{v_1,..,v_n\}$. We look for the optimizing matrix $\hat{T} \in \mathbb{R}^{(2n+1)\times(2n+1)}$ for $M(G)$ in the following form. We consider two parameters (that are taken as variables for the moment) $\delta$ and $\eta$ and write the matrix $\hat{T}$ in the following form (where $0_{k\times \ell}$ stands for a $k\times \ell$, meaning $k$ rows and $\ell$ columns, all-$0$ submatrix).
\[
\hat{T}=
\begin{bmatrix}
\frac{\delta}{|\lambda_n|}T & \frac{1}{|\lambda_n|}T & 0_{n\times 1}\\
\frac{1}{|\lambda_n|}T & 0_{n\times n} & (t-1)\sqrt{\eta}v_1\\
0_{1\times n} & (t-1)\sqrt{\eta}v_1^\top & 0
\end{bmatrix}.
\]
\medskip
\par\noindent
We can see that this matrix satisfies the conditions in Theorem~\ref{spectra} for $M(G)$, so if we can choose $\delta$ and $\eta$ so that the right hand side of the formula in Theorem~\ref{spectra} becomes equal to the left hand side of (\ref{eq:lowb}) that proves the required lower bound.

\medskip
\par\noindent
Now we define some auxiliary matrices, a $3\times 3$ one denoted $T_1$ and $n-1$ smaller, $2\times 2$ matrices that will be $T_2,\dots,T_n$. They are given as
\begin{align*}
T_1: =
\begin{bmatrix}
    (t-1)\delta & t-1 & 0\\
    t-1 & 0 & (t-1)\sqrt{\eta}\\
    0&(t-1)\sqrt{\eta}&0
\end{bmatrix}
\ {\rm and}\ \
T_i :=
\begin{bmatrix}
    \frac{\lambda_i}{|\lambda_n|}\delta& \frac{\lambda_i}{|\lambda_n|}\\
    \frac{\lambda_i}{|\lambda_n|}& 0
\end{bmatrix} \text{ for } 2\leq i \leq n.
\end{align*}

\medskip
\par\noindent
The following lemma, which was inspired by the proof of Theorem 3.1 in \cite{BKS}, will be useful.

\begin{lem}
   \[
   {\rm Sp}(\hat{T}) = \bigcup_{i=1}^n{\rm Sp}(T_i).
   \]
\end{lem}
\proof Since $T$ is a symmetric matrix we can write it in the form $A = P\Lambda P^\top$ according to Theorem \ref{symmetric}, where $P = [v_1 |...|v_n]$ and $\Lambda = {\rm diag}(\lambda_1,...,\lambda_n)$. Moreover, every $v_i$ is a unit vector and they are pairwise orthogonal. Now we use the main idea of the proof of Theorem 3.1 in \cite{BKS} for writing our matrix $\hat{T}$ in the following form.
\begin{gather*}
\begin{bmatrix}
\frac{\delta}{|\lambda_n|}T & \frac{1}{|\lambda_n|}T & 0_{n\times 1}\\
\frac{1}{|\lambda_n|}T & 0_{n\times n} & (t-1)\sqrt{\eta}v_1\\
0_{1\times n} & (t-1)\sqrt{\eta}v_1^\top & 0
\end{bmatrix}
=
\begin{bmatrix}
\frac{\delta}{|\lambda_n|}P\Lambda P^\top & \frac{1}{|\lambda_n|}P\Lambda P^\top & 0_{n\times 1}\\
\frac{1}{|\lambda_n|}P\Lambda P^\top & 0_{n\times n} & (t-1)\sqrt{\eta}v_1\\
0_{1\times n} & (t-1)\sqrt{\eta}v_1^\top & 0
\end{bmatrix}
= \\
\begin{bmatrix}
P & 0 & 0\\
0 & P &0\\
0 & 0 & 1
\end{bmatrix}
\begin{bmatrix}
\frac{\delta}{|\lambda_n|}\Lambda & \frac{1}{|\lambda_n|}\Lambda & 0_{n\times 1}\\
\frac{1}{|\lambda_n|}\Lambda & 0_{n\times n} & (t-1)\sqrt{\eta}P^\top v_1\\
0_{1\times n} & (t-1)\sqrt{\eta}v_1^\top P & 0
\end{bmatrix}
\begin{bmatrix}
P^\top & 0 & 0\\
0 & P^\top &0\\
0 & 0 & 1
\end{bmatrix}
\end{gather*}

\bigskip

\medskip
\par\noindent
Since $v_1$ is a unit vector and it is orthogonal to all other eigenvectors, $v_1^\top P = (1,0,...,0)^{\top}$ and $P^\top v_1 = (1,0,...,0)$. Hence
\[
\hat{T} = \begin{bmatrix}
P & 0 & 0\\
0 & P &0\\
0 & 0 & 1
\end{bmatrix}
\begin{bmatrix}
\frac{\delta}{|\lambda_n|}\Lambda & \frac{1}{|\lambda_n|}\Lambda & 0_{n\times 1}\\
\frac{1}{|\lambda_n|}\Lambda & 0_{n\times n} & ((t-1)\sqrt{\eta},0,..,0)\\
0_{1\times n} & ((t-1)\sqrt{\eta},0,..,0)^\top & 0
\end{bmatrix}
\begin{bmatrix}
P^\top & 0 & 0\\
0 & P^\top &0\\
0 & 0 & 1
\end{bmatrix}
\]
The spectra of
\[
\hat{\Lambda} :=
\begin{bmatrix}
\frac{\delta}{|\lambda_n|}\Lambda & \frac{1}{|\lambda_n|}\Lambda & 0_{n\times 1}\\
\frac{1}{|\lambda_n|}\Lambda & 0_{n\times n} & ((t-1)\sqrt{\eta},0,..,0)\\
0_{1\times n} & ((t-1)\sqrt{\eta},0,..,0)^\top & 0
\end{bmatrix}
\]
is identical to the spectra of $\hat{T}$ because $P$ is an orthonormal matrix.\\
So we need to prove
 \[
   {\rm Sp}(\hat{\Lambda}) = \bigcup_{i=1}^n{\rm Sp}(T_i).
   \]
We can write
$\hat{\Lambda}-\mu I$ in the following form.
\[
\left[\begin{array}{c c c c|c c c c|c}
            \frac{\delta}{|\lambda_n|}\lambda_1 - \mu & 0 & ... & 0 & \frac{1}{|\lambda_n|}\lambda_1 & 0 & ...& 0&0\\
    0 & \frac{\delta}{|\lambda_n|}\lambda_2 - \mu & ... & 0&0 & \frac{1}{|\lambda_n|}\lambda_2 &...& 0&0\\
                \vdots&\vdots&\ddots&\vdots&\vdots&\vdots&\ddots&\vdots&\vdots\\
    0 & 0 & ... & \frac{\delta}{|\lambda_n|}\lambda_n - \mu&0 & 0 & ... & \frac{1}{|\lambda_n|}\lambda_n& 0\\
    \hline
            \frac{1}{|\lambda_n|}\lambda_1 & 0 & ...& 0 & -\mu & 0 & ...& 0&(t-1)\sqrt{\eta}\\
    0 & \frac{1}{|\lambda_n|}\lambda_2  &...& 0&0 & -\mu &...& 0&0\\
                \vdots&\vdots&\ddots&\vdots&\vdots&\vdots&\ddots&\vdots&\vdots\\
    0 & 0 & ... & \frac{1}{|\lambda_n|}\lambda_n&0 & 0 & ... & -\mu& 0\\
    \hline
    0 & 0 & ... & 0 & (t-1)\sqrt{\eta} &0 & ... &0 & -\mu\\
\end{array}\right]
\]
Note that the exchanging of two rows or two columns does not change the $0$ value of a determinant and thus it leaves the roots of the characteristic polynomial unchanged. To obtain a more convenient version of $\hat{\Lambda}$ for seeing the roots of its characteristic polynomial we move its $(n+1)^{\rm th}$ row to become the second row and move the $(n+1)^{\rm th}$ column to become the second column. Similarly, we move the $(2n+1)^{\rm th}$ row to become the third row and the
$(2n+1)^{\rm th}$ column to become the third column. (We do this so that the $i^{\rm th}$ row/column becomes the $(i+2)^{\rm nd}$ when $i\notin\{1,n+1,2n+1\}$.)

\medskip
\par\noindent
The resulting matrix then looks like this:

\[
\left[\begin{array}{c c c |c c c }
            \frac{\delta}{|\lambda_n|}\lambda_1 - \mu & \frac{1}{|\lambda_n|}\lambda_1&0&0&...&0\\
            \frac{1}{|\lambda_n|}\lambda_1 & -\mu & (t-1)\sqrt{\eta}&0&...&0\\
            0&(t-1)\sqrt{\eta} & -\mu&0&...&0 \\
            \hline
            0&0&0& & & \\
            \vdots&\vdots&\vdots& & M_2 &  \\
            0&0&0& & &
\end{array}\right]
\]
where
\[
M_2 = \left[\begin{array}{c c c c | c c c c }
 \frac{\delta}{|\lambda_n|}\lambda_2 - \mu&0&...&0&\frac{1}{|\lambda_n|}\lambda_2&0&...&0\\
 0&\frac{\delta}{|\lambda_n|}\lambda_3 - \mu&...&0&0&\frac{1}{|\lambda_n|}\lambda_3&...&0\\
 \vdots &\vdots&\ddots &\vdots & \vdots &\vdots &\ddots &\vdots\\
 0&0&...&\frac{\delta}{|\lambda_n|}\lambda_n - \mu&0&0&...&\frac{1}{|\lambda_n|}\lambda_n\\
 \hline
 \frac{1}{|\lambda_n|}\lambda_2 &0& ... & 0& -\mu &0&...&0\\
 0 & \frac{1}{|\lambda_n|}\lambda_3 &...&0&0&-\mu&...&0\\
 \vdots &\vdots&\ddots&\vdots&\vdots&\vdots&\ddots&\vdots\\
 0&0&...& \frac{1}{|\lambda_n|}\lambda_n&0&0&...&-\mu
\end{array}\right]
\].

\smallskip
\par\noindent
Let
\[
M_1 := \begin{bmatrix}
            \frac{\delta}{|\lambda_n|}\lambda_1 - \mu & \frac{1}{|\lambda_n|}\lambda_1&0\\
            \frac{1}{|\lambda_n|}\lambda_1 & -\mu & (t-1)\sqrt{\eta}\\
            0&(t-1)\sqrt{\eta} & -\mu
            \end{bmatrix},
\].

\medskip
\par\noindent
Since the matrix obtained from $\hat{\Lambda}-\mu I$ with row and column order changes is block diagonal with blocks $M_1$ and $M_2$ we have ${\rm det}(\hat{\Lambda}-\mu I) = (-1)^s{\rm det}(M_1){\rm det}(M_2)$ for some appropriate positive integer $s$ that has no effect on the equation when both sides are zero. This already shows ${\rm Sp}(T_1)\subseteq {\rm Sp}(\hat{\Lambda})$. Now we manipulate further $M_2$ to get a similar product for ${\rm det}(M_2)$. To this end we move the $n^{\rm th}$ row of $M_2$ (this is the one starting with $\frac{1}{|\lambda_n|}\lambda_2$) to become the second row and the $n^{\rm th}$ column (also starting with $\frac{1}{|\lambda_n|}\lambda_2$) to become the second column. This results in the following matrix

\[
\left[\begin{array}{c c  |c c c }
                \frac{\delta}{|\lambda_n|}\lambda_2 - \mu&  \frac{1}{|\lambda_n|}\lambda_2&0&...&0\\
      \frac{1}{|\lambda_n|}\lambda_2 & -\mu&0&...&0\\
            \hline
            0&0& & & \\
            \vdots&\vdots& & M_4 &  \\
            0&0& &
\end{array}\right] ,
\]
Letting $M_3$ to be the $2\times 2$ block in the upper left corner and
\[
M_4 :=
\left[\begin{array}{c c c c | c c c c }
 \frac{\delta}{|\lambda_n|}\lambda_3 - \mu&0&...&0&\frac{1}{|\lambda_n|}\lambda_3&0&...&0\\
 0&\frac{\delta}{|\lambda_n|}\lambda_4 - \mu&...&0&0&\frac{1}{|\lambda_n|}\lambda_4&...&0\\
 \vdots &\vdots&\ddots &\vdots & \vdots &\vdots&\ddots &\vdots\\
 0&0&...&\frac{\delta}{|\lambda_n|}\lambda_n - \mu&0&0&...&\frac{1}{|\lambda_n|}\lambda_n\\
 \hline
 \frac{1}{|\lambda_n|}\lambda_3 & 0 & ... & 0& -\mu &0&...&0\\
 0 & \frac{1}{|\lambda_n|}\lambda_4 &...&0&0&-\mu&...&0\\
 \vdots &\vdots&\ddots&\vdots&\vdots&\vdots&\ddots&\vdots\\
 0&0&...& \frac{1}{|\lambda_n|}\lambda_n&0&0&...&-\mu
\end{array}\right].
\]
we get
${\rm det}(M_2) = (-1)^{s'}{\rm det}(M_3){\rm det}(M_4)$ (where again, $s'$ is irrelevant when both sides are equal to $0$) and since ${\rm det}(M_2)$ is just the characteristic polynomial of $T_2$, this shows ${\rm Sp}(T_2)\subseteq {\rm Sp}(\hat{\Lambda})$. Continuing in a similar fashion until the last remaining block has also size $2\times 2$ we obtain $\cup_{i=1}^n {\rm Sp}(T_i)={\rm Sp}(\hat{\Lambda})$, completing the proof of the lemma. \qed    \\

\medskip
\par\noindent
\vspace{0.5 cm}
Observe that $T_1$ and $T_n$ can be written in the simpler forms $T_1 = (t-1)T_1^*$, where
\[
T_1^* =
\begin{bmatrix}
    \delta & 1 & 0\\
    1 & 0 & \sqrt{\eta}\\
    0&\sqrt{\eta}&0
\end{bmatrix},\]
and
\[
T_n =
\begin{bmatrix}
    -\delta& -1\\
    -1& 0
\end{bmatrix}.
\]

Our plan is to show that there exist some $\delta,\eta,\gamma \in \mathbb{R}_+$ such that $\lambda_{{\rm max}}(\hat{T})= \gamma(t-1)$ and $\lambda_{{\rm min}}(\hat{T}) = -\gamma\frac{t-1}{m(t)-1}$. If these two conditions are satisfied, then $m(t)=1+\frac{\lambda_{{\rm max}}(\hat{T})}{|\lambda_{{\rm min}}(\hat{T})|}$.

\smallskip
\par\noindent
In order to determine the appropriate $\delta,\eta,\gamma$, we need to calculate the characteristic polynomials of these matrices, however, as we will see that is not quite necessary. It will be sufficient if we determine the characteristic polynomials of $T_1^*$ and $T_n$,

\smallskip
\par\noindent
Let
$p_1^*(\mu)={\rm det}(T_1^*-\mu I)$, $p_n(\mu):={\rm det}(T_n-\mu I)$. We will also use $p_1(\mu) := {\rm det}\left(T_1 - \mu I\right)$ as an auxiliary polynomial.
The roots of $p_1^*(\mu)$ are the solutions of $\mu^3 -\delta\mu^2 -(\eta +1)\mu + \eta\delta=0$ (the left hand side being just $-p_1^*(\mu)$). These values
and the solutions of $p_n(\mu) = \mu^2 + \delta\mu -1= 0$ will be the eigenvalues of $\hat{T}$ that we need.\\

\medskip
\par\noindent
We want the minimum and maximum eigenvalues of $\hat{T}$ to be
$\lambda_{{\rm min}}(\hat{T}) = -\gamma\frac{t-1}{m(t)-1}$ and
$\lambda_{{\rm max}}(\hat{T})= \gamma(t-1)$, respectively. This would be enough for completing the proof, since by Theorem~\ref{spectra} it would give us $\thet(G)\ge m(t)$, while $\thet(G)\le m(t)$ we have already shown.

\medskip
\par\noindent
We will see that the minimum and maximum eigenvalues will be as chosen above if the following equations are satisfied.
\begin{enumerate}
    \item $p_n\left(-\gamma\frac{t-1}{m(t)-1}\right) = 0$;
    \item $p_1\left(\gamma(t-1)\right) = 0$ $\Leftrightarrow$ $p_1^*(\gamma) = 0$;
    \item $p_1\left(-\gamma\frac{t-1}{m(t)-1}\right) = 0$ $\Leftrightarrow$ $p_1^*\left(-\gamma\frac{1}{m(t)-1}\right)=0$.
\end{enumerate}

\medskip
\par\noindent
This means that we fix $-\gamma\frac{t-1}{m(t)-1}$ to be an eigenvalue of both $T_1$ and $T_n$, while $\gamma(t-1)$ is also an eigenvalue of $T_1$. On the way it will also be made sure that $\gamma>0$ and that no matrix $T_i$ in our collection has a smaller eigenvalue than $-\gamma\frac{t-1}{m(t)-1}$, that is the absolute value of $\lambda_{\min}$ is indeed what this value provides. By this the maximality of the other eigenvalue $\gamma(t-1)$ will already be guaranteed, since otherwise we would get $1+\frac{\lambda_{\max}(\hat{T})}{|\lambda_{\min}(\hat{T})|}$ larger than $m(t)$, which is impossible by $m(t)\ge\thet(G)$ and Theorem~\ref{spectra}.

\medskip
\par\noindent
\smallskip
\par\noindent
Assume that $\mu_1$ and $\mu_2$ are the two roots of $p_n(\mu)$, therefore $\mu_1\mu_2 = -1$ and $-(\mu_1 + \mu_2)=\delta$. Thus fixing $\mu_2$ to be equal to $-\gamma\frac{t-1}{m(t)-1}$ gives $\mu_1 = \frac{1}{\gamma\frac{t-1}{m(t)-1}}$ from which we can calculate $\delta$:
\[
\delta = \gamma\frac{t-1}{m(t)-1} - \frac{1}{\gamma\frac{t-1}{m(t)-1}}.
\]
Now we have to make sure that the choice of $\delta$ still lets the value
$\mu_2=-\gamma\frac{t-1}{m(t)-1}$ being the smallest eigenvalue of $\hat{T}$, that is, that no eigenvalue of $T_i$ for $2\leq i \leq n-1$ becomes smaller than $\mu_2$. So we want to guarantee ${\rm Sp}(T_i) \subseteq \left[-\gamma\frac{t-1}{m(t)-1},\infty \right[$.\\

    \medskip
    \par\noindent
    Recall that $\lambda_1\ge\dots\ge\lambda_n$ are the eigenvalues of $T$.

    \smallskip
    \par\noindent
    \textbf{Case 1:} If $\lambda_i \leq 0$, then using that
    \begin{align*}\label{transform}
T_i=
      \begin{bmatrix}
    \frac{\lambda_i}{|\lambda_n|}\delta& \frac{\lambda_i}{|\lambda_n|}\\
    \frac{\lambda_i}{|\lambda_n|}& 0
\end{bmatrix}=
-\frac{\lambda_i}{|\lambda_n|}
\begin{bmatrix}
    -\delta& -1\\
    -1& 0
\end{bmatrix}=\frac{\lambda_i}{\lambda_n}T_n
    \end{align*}
    and that $\lambda_n\le \lambda_i\le 0$
     we can write $\lambda_{{\rm min}}(T_i)=\frac{\lambda_i}{\lambda_n}\lambda_{\rm min}(T_n)=\gamma\frac{\lambda_i}{|\lambda_n|}\frac{t-1}{m(t)-1} \geq \gamma\frac{\lambda_n}{|\lambda_n|}\frac{t-1}{m(t)-1} = -\gamma\frac{t-1}{m(t)-1}$.\\

\textbf{Case 2:} Assume that $\lambda_i > 0$.
Then we have $\frac{\lambda_i}{\lambda_n}<0$ by $\lambda_n<0$ (for the latter cf. the discussion after Theorem~\ref{spectra}),
so we can write $\lambda_{{\rm min}}(T_i)=-\frac{\lambda_i}{\lambda_n}\lambda_{\rm max}(T_n)=-\frac{\lambda_i}{|\lambda_n|}\frac{1}{\gamma\frac{t-1}{m(t)-1}} \geq -\frac{\lambda_1}{|\lambda_n|}\frac{1}{\gamma\frac{t-1}{m(t)-1}} = -(t-1)\frac{1}{\gamma\frac{t-1}{m(t)-1}}$. Thus the required inequality
$\lambda_{\rm min}(T_i)\ge \lambda_{\rm min}(T_n)$ follows if the following inequality holds.
$$ -(t-1)\frac{1}{\gamma\frac{t-1}{m(t)-1}}\ge -\gamma\frac{t-1}{m(t)-1}$$
\smallskip
\par\noindent
The latter inequality will follow from $(m(t)-1)^2\leq\gamma^2(t-1)$ if $\gamma>0$ also holds.
The proof of these two inequalities we postpone to the final part of the proof, see Lemma~\ref{lem:maradt} and its proof in Appendix 3.

\bigskip
\par\noindent
We continue by determining the value of $\eta$ and $\gamma$.
Let $\mu_1,\mu_2,\mu_3$ be the solutions of the $\mu^3 -\delta\mu^2 -(\eta +1)\mu + \eta\delta=0$ equation we obtained from $p_1^*=0$. We know that $\mu_1 = \gamma$ and $\mu_3 = -\gamma\frac{1}{m(t)-1}$. Thus by Viete's formulas describing the relationship of the roots and coefficients of a polynomial, we get the following system of equations.
\begin{equation}
\centering
\left\{\begin{split}
\mu_1\mu_2\mu_3 &= -\eta\delta \\
\mu_1 + \mu_2 + \mu_3 &= \delta \\
\mu_1\mu_2 + \mu_1\mu_3 + \mu_2\mu_3 &= -(\eta + 1).\
\end{split}\right.
 \end{equation}
 From the first equation, we have $\mu_2 = \frac{1}{\gamma^2}(m(t)-1)\eta\delta$, which gives $0 \leq \mu_2$ if we also have $0 \leq \eta$ and $0 \leq \delta$ both of which will be shown at the end of the proof. (See Lemma~\ref{lem:maradt} and its proof in Appendix 3.) Taking this for granted now, we know that $\mu_2$ is surely not the minimal eigenvalue of $\hat{T}$. Now we substitute the above value of $\mu_2$ into the second equation and solve it for $\eta$. That gives us
 \[
 \eta = \frac{1}{\delta(m(t)-1)}\left(\gamma^2\delta + \gamma^3\frac{1}{m(t)-1}-\gamma^3\right).
 \]
(For the last equation to make sense we will actually need the strict inequality $\delta>0$. To see that this will be satisfied we refer again to Lemma~\ref{lem:maradt} and its proof in Appendix 3.)

\medskip
\par\noindent
Substituting the above expressions for $\mu_2$ and $\eta$ into the third equation we obtain
\begin{gather*}
\gamma\delta + \gamma^2\frac{1}{m(t)-1}-\gamma^2-\gamma^2\frac{1}{m(t)-1}-\frac{1}{m(t)-1}\gamma\delta-\gamma^2\frac{1}{(m(t)-1)^2} + \gamma^2\frac{1}{m(t)-1}=\\
-\ \frac{1}{\delta(m(t)-1)}\left(\gamma^2\delta + \gamma^3\frac{1}{m(t)-1}-\gamma^3\right) -1.
\end{gather*}
We have already seen that $\delta = \gamma\frac{t-1}{m(t)-1} - \frac{1}{\gamma\frac{t-1}{m(t)-1}}$, so this equation depends only on $\gamma, t$ and $m(t)$. (Although $m(t)$ is just a function of $t$, here we handle it as a parameter itself.)

The last equation can be transformed into one having a polynomial with variable $\gamma$ and coefficients depending solely on $t$ and $m(t)$ on the left had side and $0$ on the right hand side. The details of this calculation can be found in Appendix 4. To give the resulting equation we use the notation $v:=t-1$ and $w:=\frac{1}{m(t)-1}$. What we obtain is $a\gamma^4 + b\gamma^2 + c = 0$, where
\begin{align*}
    a =& -v^4w^5 + v^4w^4 - v^3w^5 + 2v^3w^4 - v^3w^3 + v^2w^4 - v^2w^3,\\
    b =& v^3w^3 + 2v^2w^3 - 2v^2w^2 + vw^3 - 2vw^2 + vw,\\
    c =& - vw - w + 1.
\end{align*}

Let us substitute $\hat{\gamma}:=\gamma^2$ to make our equation quadratic: $$a\hat{\gamma}^2 + b\hat{\gamma} + c = 0.$$
The discriminant $D(t,m(t)) =b^2-4ac$ of this quadratic equation turns out to be equal to $0$.
In fact, the discriminant is equal to
\begin{gather*}
    (v^3w^3 + 2v^2w^3 - 2v^2w^2 + vw^3 - 2vw^2 + vw)^2\\
    -4( -v^4w^5 + v^4w^4 - v^3w^5 + 2v^3w^4 - v^3w^3 + v^2w^4 - v^2w^3)(- vw - w + 1) = .
\end{gather*}
\[
v^6w^6 - 2v^4w^6 + 4v^4w^5 - 2v^4w^4 + 4v^3w^5 - 4v^3w^4 + v^2w^6 - 2v^2w^4 + v^2w^2 =
\]
\[
v^2 w^2 (v w + w - 1) (v^3 w^3 - v^2 w^3 + v^2 w^2 - v w^3 + 2 v w^2 - v w + w^3 + w^2 - w - 1).
\]

\medskip
\par\noindent
Hence $D(t,m(t))=0$ will follow if we have $v^3 w^3 - v^2 w^3 + v^2 w^2 - v w^3 + 2 v w^2 - v w + w^3 + w^2 - w - 1 = 0$.
\medskip
\par\noindent
Substituting back $t-1$ for $v$ and $\frac{1}{m(t)-1}$  for $w$ in the last equation, after some manipulation we get the equation
$m(t)^3 + (t-3)m(t)^2 + (3-2t-t^2)m(t) + (-t^3 +5t^2 -3t-1) = 0$ which is the same as equation (\ref{eq:cubic}) in the first half of the proof if we exchange $m(t)$ to $\hat{t}$. Since $m(t)$ was defined as its solution for $\hat{t}$, we know that this equation holds, therefore we indeed have $D(t,m(t))=0$.
\medskip
\par\noindent
Since $D(t,m(t))$ was the discriminant of the quadratic equation $a\hat{\gamma}^2 + b\hat{\gamma} + c = 0$, we now know that this equation has a unique real solution which is the following.
\begin{gather*}
\hat{\gamma} = -\frac{b}{2a} = -\frac{v^3w^3 + 2v^2w^3 - 2v^2w^2 + vw^3 - 2vw^2 + vw}{2(-v^4w^5 + v^4w^4 - v^3w^5 + 2v^3w^4 - v^3w^3 + v^2w^4 - v^2w^3)}=\\
-\frac{(m(t) - 1)^2 (t - m(t) + 1)^2}{2 (t - 1) t (m(t) - 2) (t - m(t))}.
\end{gather*}

\medskip
\par\noindent
Since $\hat{\gamma}$ was defined as $\gamma^2$, we have to check that we obtained a non-negative value for it. Since $t\ge 2$ and $m(t)\ge t$, this will follow if we show that the last inequality is a strict one, that is, $m(t)>t$. To exclude $m(t)=t$ it is enough to observe that the value of the left hand side of our cubic equation (\ref{eq:cubic}) for $m(t)=t$ is constant $-1$. Therefore there does not exist any $t$ such that $t=m(t)$. So we obtained a meaningful solution for $\hat{\gamma}$.
In fact, we also see from the formula for it, that it is strictly positive.
Since we have the right to choose the positive root of $\hat{\gamma}$ to be $\gamma$, the inequality $\gamma>0$ we needed earlier is also settled.

\bigskip

\bigskip
\par\noindent
We still owe the proof of some inequalities we have taken for granted on the way postponing their proof to the end. These are the following.

\begin{lem}\label{lem:maradt}
    The following inequalities are true:
    \begin{enumerate}
    \item $0< \eta$;
    \item $\frac{(m(t)-1)^2}{t-1} \leq \gamma^2$;
    \item $0< \delta$.
    \end{enumerate}
\end{lem}

The proof of this Lemma is given in Appendix 3.

\bigskip
\par\noindent
By the above we have completed the proof of Theorem~\ref{main}. $\Box$

\section{Further variants}\label{sect:egyebek}

\subsection{Sperner capacity of a Mycielskian}

\medskip
\par\noindent
Sperner capacity, first defined in \cite{GKV1}, generalizes the notion of Shannon capacity to directed graphs. It also has an operational meaning in terms of information transmission with zero-error over a compound channel, see \cite{NayRose} for details.

\medskip
\par\noindent
To give the definition of Sperner capacity we have to generalize the OR-product to directed graphs (often called simply digraphs). This is rather straightforward.

\medskip
\par\noindent
\begin{defi} \label{defi:ORdir}
For two digraphs $F$ and $G$ their OR-product $F\cdot G$ is defined by
$$V(F\cdot G)=V(F)\times V(G)$$ and
$$E(F\cdot G)=\{((f,g),(f',g')): f,f'\in V(F), g,g'\in V(G),$$$$(f,f')\in E(F)\ {\rm or}\ (g,g')\in E(G)\}.$$
The $t^{\rm th}$ OR-power $G^t$ of a digraph $G$ is meant to be the $t$-fold OR-product of $G$ with itself.
\end{defi}

\medskip
\par\noindent
It can be seen that the difference between Definition~\ref{defi:ORprod} and Definition~\ref{defi:ORdir} is only that unordered pairs of vertices are exchanged to ordered ones. Note that even if both $F$ and $G$ are oriented graphs (meaning digraphs that have no cycles of length $2$, that is, for no edge its reversed edge is present), their product may contain $2$-cycles. Let $T_n$ denote the transitive tournament on $n$ vertices, that is the oriented complete graph, whose vertices can be labeled by $1,2,\dots,n$ so that $(i,j)\in E(T_n)\Leftrightarrow i<j$. The subdigraph of a digraph $D$ induced by $U\subseteq V(D)$ will be denoted $D[U]$. A subdigraph isomorphic to the transitive tournament $T_m$ is called a transitive clique of size $m$.

\medskip
\par\noindent
\begin{defi}\label{defi:symmtr}
  For a digraph $D$ its {\em symmetric clique number} is defined as
  $$\omega_s(D):=\max\{|Q|: Q\subseteq V(D), \forall u,v\in Q, (u,v),(v,u)\in E(D)\}.$$
  The {\em transitive clique number} of $D$ is defined as
  $$\omega_{\rm tr}(D):=\max\{|Q|: Q\subseteq V(D), T_{|Q|}\subseteq D[Q]\}.$$
\end{defi}

\medskip
\par\noindent
Note that we have $\omega_s(D)\le\omega_{\rm tr}(D)$ for any digraph $D$.
The straighforward generalization of Shannon OR-capacity to digraphs is obtained if we consider undirected graphs as directed ones where all edges are present in both directions, thus exchange $\omega(G)$ to $\omega_s(G)$ in Definition~\ref{defi:orcap}. This is how Sperner capacity is defined in \cite{GKV1}. It is not hard to see, however, that the so defined limit remains the same if we change $\omega_s(D^t)$ to $\omega_{\rm tr}(D^t)$ in the definition. This alternative definition of Sperner capacity appears already, e.g., in \cite{Alondir}. For us this way of defining it is also more convenient.

\medskip
\par\noindent
\begin{defi}\label{defi:spcap}
  The (non-logarithmic) Sperner capacity of the digraph $D$ is defined by
  $$C_{\rm Sp}(D):=\lim_{t\to\infty}\sqrt[t]{\omega_{\rm tr}(D^t)}.$$
\end{defi}

\medskip
\par\noindent
Like in case of Shannon capacity, it is straightforward to see that $\sqrt[k]{\omega_{\rm tr}(D^k)}\le C_{\rm Sp}(D)$ holds for any fixed positive integer $k$ and any digraph $D$, from which the existence of the limit follows again by Fekete's Lemma.
We say that the Sperner capacity of digraph $D$ is attained in finite length if there exists some positive integer $k$ for which $C_{\rm Sp}(D)=\sqrt[k]{\omega_{\rm tr}(D^k)}$.

\medskip
\par\noindent
Now we define the Mycielskian of a digraph.

\medskip
\par\noindent
\begin{defi}\label{defi:dirMyc}
Let $D$ be a digraph.
Its Mycielskian $M(D)$ is defined on the vertices $$V(M(D))=V(D)\times\{0,1\}\cup \{z_{M(D)}\}$$
with edge set $$E(M(D))=\{((v,0)(w,i)): (v,w)\in E(G), i\in\{0,1\}\}\cup$$
$$\{((v,1)(w,0)): (v,w)\in E(G)\}\}\cup\{(z_{M(D)},(v,1)): v\in V(D)\}.$$
\end{defi}

\medskip
\par\noindent
Thus the edges between the ``levels'' $V(D)\times\{0\}$ and $V(D)\times\{1\}$ inherit the orientation of the corresponding edges in $E(D)$ while all edges between the vertices in $V(D)\times\{1\}$ and the vertex $z_{M(D)}$ are oriented outward from $z_{M(D)}$. (Note that although orienting all the edges between $z_{M(D)}$ and its neighbors towards $z_{M(D)}$ would be equally natural, so this choice is just a matter of convention, it will not alter essentially anything in what follows. The reason is the simple observation that if we reverse all the edges of a directed graph then its Sperner capacity does not change.)

With this definition the following extension of
Theorem~\ref{thm:fintcap} to Sperner capacity is almost straightforward.

\medskip
\par\noindent
\begin{thm}\label{thm:fintsp}
  If $D$ is a digraph that attains its Sperner capacity in finite length, then
  $$C_{\rm Sp}(M(D))>C_{\rm Sp}(D).$$
\end{thm}

\smallskip
\par\noindent
The directed generalization of Lemma~\ref{lem:felemel} is the following.

\medskip
\par\noindent
\begin{lem}~\label{lem:dirfelemel}
There is a transitive clique of size $n^n$ in $[M(T_n)\setminus\{z_{M(T_n)}\}]^n$ every vertex of which has a coordinate belonging to the set $V(T_n)\times\{1\}$.
\end{lem}

\medskip
\par\noindent
\proof
Let the vertices of $T_n$ be $0,1,\dots,n-1$ so that $(i,j)\in E(T_n)$ if $i<j$.
Define the set $\hat{B}_n$ the same way as it is done in the proof of Lemma~\ref{lem:felemel}. We claim that the subdigraph induced by $\hat{B}_n$ contains a transitive clique of size $n^n$.
For a vertex $({\mbf x},{\mbf i}):=(x_1,i_1)(x_2,i_2)\dots (x_n,i_n)$ with $i_j\in\{0,1\}\ \forall j$ let $s({\mbf x},{\mbf i})=\sum_{j=1}^n x_j$. If for $({\mbf x},{\mbf i}), ({\mbf x'},{\mbf i'})\in \hat{B}_n$ we have $s({\mbf x},{\mbf i})<s({\mbf x'},{\mbf i'})$ then we have the directed edge
$(({\mbf x},{\mbf i})({\mbf x'},{\mbf i'}))$ present in $(M(T_n))^n[\hat{B}_n]$. This can be seen as follows. Either
$s({\mbf x},{\mbf i})\equiv s({\mbf x'},{\mbf i'})\pmod n$ or not. In the latter case there must exist a coordinate $j$ where $x_j<x'_j$ and at most one of $i_j$ and $i'_j$ is equal to $1$, therefore $((x_j,i_j)(x'_j,i'_j))\in E(M(T_n))$. In the former case there must be at least two coordinates $j$ where $x_j<x'_j$, since $s({\mbf x'},{\mbf i'})-s({\mbf x},{\mbf i})\ge n$ in this case, while $x'_j-x_j$ can be at most $n-1$. But one of these two coordinates will also satisfy $i_j=i'_j=0$ and so $((x_j,i_j)(x'_j,i'_j))\in E(M(T_n))$ follows again.

\medskip
\par\noindent
Taking into account the above, it is enough to prove that all vertices $({\mbf x},{\mbf i})\in \hat{B}_n$ for which $s({\mbf x},{\mbf i})=k$ for an arbitrary constant $k$ also induce a subdigraph of
$[M(T_n)]^n$ which contains a transitive tournament on all these vertices. Note that $s({\mbf x},{\mbf i})=k=s({\mbf x'},{\mbf i'})$ for vertices in $\hat{B}_n$ implies that ${\mbf i}={\mbf i'}$, thus there is a unique coordinate $r$ for which $i_r=i'_r=1$ while $i_j=i'_j=0$ for all $j\neq r$. The $s({\mbf x},{\mbf i})=s({\mbf x'},{\mbf i'})$ relation also implies that ${\mbf x}$ and ${\mbf x'}$ differs at least at two coordinates, so there is some $j\neq r$ where $x_j\neq x'_j$. Now we define an ordering of the vertices just considered for which an earlier vertex in this order always sends an edge to the ones coming later in this order, thus the transitive tournament we seek for is found. Put $({\mbf x},{\mbf i})$ in front of $({\mbf x'},{\mbf i'})$ if for the first $j\neq r$ where they differ, $x_j<x'_j$. Since for such a $j$ we have $i_j=i'_j=0$ we will have $(({\mbf x},{\mbf i}),({\mbf x'},{\mbf i'}))\in E([M(T_n)]^n)$ as we need. Note that by this rule we also do not put any three vertices into cyclic order, that is, if $({\mbf x},{\mbf i})$ precedes $({\mbf x'},{\mbf i'})$ and $({\mbf x'},{\mbf i'})$ precedes $({\mbf x''},{\mbf i''})$ then $({\mbf x},{\mbf i})$ also precedes $({\mbf x''},{\mbf i''})$.
The definition of $\hat{B}_n$ also ensures that each vertex $({\mbf x},{\mbf i})$ in it has a coordinate $j$ with $i_j=1$. This completes the proof.
\qed

\medskip
\par\noindent
The above lemma implies the following generalization of Corollary~\ref{lem:Kn}, which is a special case of Theorem~\ref{thm:fintsp}.
\medskip
\par\noindent
\begin{cor}\label{lem:Tn}
  For the transitive tournament $T_n$ on $n$ vertices we have
  $$C_{\rm Sp}(M(T_n))\ge \sqrt[n]{n^n+1}>n.$$
 \end{cor}

\proof
Consider the transitive clique of size $n^n$ in $[M(T_n)\setminus\{z_{M(T_n)}\}]^n$ whose existence is shown in Lemma~\ref{lem:dirfelemel}. Since all the sequences that are vertices of this clique contain a coordinate from $V(T_n)\times \{1\}=N_+(z_{M(T_n)})$ (where $N_+(v)$ stands for the outneighborhood of vertex $v$), the all-$z_{M(T_n)}$ sequence of length $n$ sends an edge to all vertices of this clique in $[M(T_n)]^n$. So adding this sequence to our transitive clique of size $n^n$ we obtain a transitive clique with $n^n+1$ vertices, thus showing $\omega_{\rm tr}([M(T_n)]^n)\ge n^n+1$. This implies $C_{\rm Sp}(M(T_n))\ge\sqrt[n]{n^n+1}$ as stated.
\qed

\medskip
\par\noindent
\begin{remark}\label{rem:spC5}
We remark that for $n=2$ $M(T_2)$ is a $5$-cycle with exactly one vertex having outdegree equal to $1$ while all the other four vertices have outdegree $2$ or $0$. This is the unique orientation of the $5$-cycle $C_5$ having Sperner capacity strictly larger than $2=C_{\rm Sp}(T_2)$ and, in fact, its Sperner capacity is equal to $\sqrt{5}=C_{\rm OR}(C_5)$. The latter is shown in \cite{GGKS}, cf. also \cite{SaSi} for a generalization. The fact, that all other orientations have Sperner capacity $2$ is already remarked in \cite{GGKS}, for the proof of a more general theorem see \cite{KPS}.
\end{remark}

\medskip
\par\noindent
Note that Lemma~\ref{lem:benne} remains true for directed graphs, so
$$M(D^t)\subseteq [M(D)]^t$$
for all directed graphs $D$.
(The proof is literally the same, so we omit the details.)

\medskip
\par\noindent
From here it is easy to complete the proof of Theorem~\ref{thm:fintsp} analogously to that of Theorem~\ref{thm:fintcap}. (It is essentially the same argument as that proving Theorem~\ref{thm:fintcap} that we repeat only for the sake of clarity and completeness.)

\medskip
\par\noindent
{\em Proof of Theorem~\ref{thm:fintsp}.}
Let $D$ be a digraph with $C_{\rm Sp}(G)=\sqrt[k]{\omega_{\rm tr}(D^k)}=\sqrt[k]{N}$ where $k\ge 1$ is a positive integer and $N=\omega_{\rm tr}(D^k)$.
Since $T_N\subseteq D^k$, we also have $M(T_N)\subseteq M(D^k)\subseteq [M(D)]^k,$ where the last relation is by the digraph version of Lemma~\ref{lem:benne}.
This further implies $$[M(T_N)]^N\subseteq [M(D^k)]^N\subseteq [M(D)]^{kN},$$
so by the proof of Corollary~\ref{lem:Tn} we can write
$$N^N+1\le \omega_{\rm tr}([M(T_N)]^N)\le\omega_{\rm tr}([M(D^k)]^N)\le\omega_{\rm tr}([M(D)]^{kN}).$$
The latter implies
$$C_{\rm Sp}(M(D))\ge \sqrt[kN]{\omega_{\rm tr}([M(D)]^{kN})}\ge\sqrt[kN]{N^N+1}>\sqrt[k]{N}=C_{\rm Sp}(D)$$
completing the proof.
\qed

\subsection{Generalized Mycielskian}

\bigskip
\par\noindent
The generalized Mycielski construction was first considered by Stiebitz \cite{Stieb}, cf. also \cite{GyJS, CTcone} and the book \cite{Mat}. Below we give the definition to directed graphs that simplifies for undirected graphs if we disregard the orientation of the edges.

\medskip
\par\noindent
\begin{defi}\label{defi:dirgenMyc}
Let $D$ be a digraph.
Its $r$-level generalized Mycielskian $M_r(D)$ is defined on the vertices $$V(M(D))=V(D)\times\{0,1,\dots,r-1\}\cup \{z_{M_r(D)}\}$$
with edge set $$E(M(D))=\{((v,0)(w,0)): (v,w)\in E(D)\}\cup$$
$$\{((v,i)(w,i+1)): (v,w)\in E(D), i\in\{0,1,\dots r-2\}\}\cup$$
$$\{((v,i+1)(w,i)): (v,w)\in E(D), i\in\{0,1,\dots r-2\}\}\cup   \{(z_{M_r(D)},(v,r-1)): v\in V(D)\}.$$
\end{defi}

\medskip
\par\noindent
Note that we get back the Mycielskian $M(D)$ for $r=2$, while $M_1(D)$ can be considered as a copy of $D$ with all its vertices connected to a new vertex (and all these edges oriented towards the vertex in $V(D)$ in the directed case).
An interesting fact about the generalized Mycielskian $M_r(G)$ of an undirected graph is that although for several graphs $G$ we have $\chi(M_r(G))=\chi(G)+1$ for all $r\ge 1$, there are also examples of graphs $G$ for which $\chi(M_r(G))=\chi(G)$ for large enough $r>2$. An example of the latter is the complementary graph of the $7$-cycle already with $r=3$, see \cite{CTcone}. (Cf. also \cite{PZ} or the conference paper \cite{GST} written to popularize a special case of the more general result in the former paper.)
Note that $M_r(K_2)\cong C_{2r+1}$. This immediatley implies that the chromatic number does increase for every generalized Mycielskian of a bipartite graph (since it will contain an odd cycle) and Stiebitz proved that the generalized Mycielskian of any odd cycle will have chromatic number $4$, that is $1$ larger than the chromatic number of the odd cycle itself. (All this is related to the topological lower bound on the chromatic number introduced by Lov\'asz \cite{LLKn}, cf. \cite{Mat} for more details.)

\medskip
\par\noindent
In the light of the above it is quite natural to ask what we can say about the effect of the generalized Mycielski construction on Shannon OR-capacity and Sperner capacity. Below we give some preliminary observations on this likely to be difficult problem. Probably the very first question that one should be able to answer is whether $C_{\rm OR}(C_{2r+1})=C_{\rm OR}(M_r(K_2))>C_{\rm OR}(K_2)=2$ for every positive integer $r$. For this question an affirmative answer is already known, but even this was far from trivial: it needed the ingenious construction of Bohman and Holzman that is the main result of \cite{BH} (presented in the usual complementary language). Their construction generalizes the one by Shannon~\cite{Sha56} proving $C_{\rm OR}(C_5)\ge\sqrt{5}$: for large enough $t$ (in particular, for $t=2^{2^{r-1}}$) they construct a clique of size $2^t$ in $[M_r(K_2)\setminus\{z_{M_r(K_2)}\}]^t$ every vertex of which contains a coordinate from $V(K_2)\times\{r-1\}$. Therefore the all-$z_{M_r(K_2)}$ sequence of length $t$ can be added to it to obtain a clique of size $2^t+1$ in the $t^{\rm th}$ OR-power, thereby proving $C_{\rm OR}(M_r(K_2))>2$. Note that our proof of Corollary~\ref{lem:Kn} uses a similar idea: we presented a clique of size $n^t$ in
$[M(K_n)\setminus\{z_{M(K_n)}\}]^t$ for large enough $t$ (in our case $t=n$ was sufficiently large) in such a way that all the vertices in it contained a coordinate from $V(K_n)\times\{1\}$ resulting in the same conclusion, that the all-$z_{M(K_n)}$ sequence of length $t$ could be added to obtain a clique of size $n^t+1$ thereby proving $C_{\rm OR}(M(K_n))>n=C_{\rm OR}(K_n)$. Our first observation is that a construction of this type cannot exist for $M_r(K_n)$ if both $r>2$ and $n>2$.

\medskip
\par\noindent
\begin{prop}\label{prop:none}
For $n,r\ge 3$ there exists no clique of size $n^t$ in $[M_r(K_n)\setminus\{z_{M_r(K_n)}\}]^t$ for any positive integer $t$ such that each of its vertices would contain a coordinate from $V(K_n)\times\{r-1\}$.
\end{prop}

\proof
Assume for contradiction that for some $n,r\ge 3$ and positive integer $t$ there exists such a clique $Q_t$. For a vertex $({\mbf x},{\mbf i}_{\msbf x})\in [M_r(K_n)\setminus\{z_{M_r(K_n)}\}]^t$, ${\mbf x}\in [V(K_n)]^t, {\mbf i}_{\msbf x}\in \{0,1,\dots, r-1\}^t$ let us call ${\mbf x}$ the projection of this vertex to $[V(K_n)]^t$. Observe that if we take the projection of all vertices in $Q_t$ to $[V(K_n)]^t$, then they must all be different, therefore give us every sequence of $[V(K_n)]^t$ exactly once. Consider some $({\mbf x},{\mbf i_{\msbf x}})\in Q_t$ and assume that $j$ is a coordinate where $i_{{\msbf x},j}=r-1$. Since $n\ge 3$ there exists two sequences $({\mbf y}, {\mbf i}_{\msbf y}), ({\mbf y'}, {\mbf i}_{\msbf y'})\in Q_t$ whose projections to $[V(K_n)]^t$ differ from the projection of $({\mbf x}, {\mbf i}_{\msbf x})$ (and thus also from the projection of each other) only in the $j^{\rm th}$ coordinate. This means that for being adjacent in $[M_r(K_n)]^t$ to $({\mbf x},{\mbf i}_{\msbf x})$ the $j^{\rm th}$ coordinate of both ${\mbf i}_{\msbf y}$ and ${\mbf i}_{\msbf y'}$ must be equal to $r-2$. But then $r-1\ge 2$ implies that $({\mbf y}, {\mbf i}_{\msbf y})$ and $({\mbf y'}, {\mbf i}_{\msbf y'})$ cannot be adjacent to each other. This contradicts the assumption that $Q_t$ is a clique.
\qed

\medskip
\par\noindent
It is straightforward from the definition of Sperner capacity that if $\vec{G}$ is a directed version of graph $G$, then $C_{\rm Sp}(\vec{G})\le C_{\rm OR}(G)$. Therefore it follows from Proposition~\ref{prop:none} that Lemma~\ref{lem:dirfelemel} cannot be generalized to $T_n$ in place of $K_n$ if $n,r\ge 3$. (Note that this does not mean that Corollary~\ref{lem:Tn} itself cannot be generalized in a different way.) This leaves the question of whether such a generalization is possible if $r\ge 3$, but $n=2$, the case for which the undirected case is settled by the result of Bohman and Holzman \cite{BH}. The situation here is rather puzzling. It is proven in \cite{KPS} that all oriented versions except possibly one of an odd cycle ({\rm oriented} meaning again that every edge is present in exactly one direction) has Sperner capacity equal to $2<C_{\rm OR}(C_{2r+1})$. The exceptional orientation, that is called {\em alternating} in \cite{KPS} is just the one that is isomorphic to $M_r(T_2)$. (The name ``alternating'' refers to the property of this kind of oriented odd cycles that their pairs of edges meeting at a common vertex are oriented in the opposite direction at as many vertices as possible, meaning at all vertices but one.) It is already shown in \cite{KPS} that for the alternating $7$-cycle $M_3(T_2)$ we do have $C_{\rm Sp}(M_3(T_2))>2$ by just adapting the Bohman-Holzman construction to the oriented case. (It also works for the alternating $C_5$ but that belongs to the case $r=2$ which we already discussed in Remark~\ref{rem:spC5}.) The puzzling fact is that (as also discussed in \cite{KPS}) the same method fails for $C_9$, that is, for $r=4$ and therefore for all $r\ge 4$. (The latter is a consequence of the fact that any alternating odd cycle maps homomorphically to any shorter alternating odd cycle.) Recently we have put some more effort into trying to adapt the ideas of Bohman and Holzman for the $C_9$ case less directly but we still failed. But we also do not have a negative result analogous to Proposition~\ref{prop:none} for the oriented case when $n=2$ and $r\ge 4$. Let us emphasize that even if we had that would not imply anything about the question whether $C_{\rm Sp}(M_r(T_2))$ can be larger than $2$ for $r\ge 4$. These questions remain the subject of further research.

\medskip
\par\noindent
Concerning the generalized Mycielski construction we also mention that Tardif \cite{CTcone} proved a generalization of the Larsen-Propp-Ullman theorem of \cite{LPU} and showed that the fractional chromatic number of $M_r(G)$ is also determined by the value of the fractional chromatic number of $G$. Whether an analogous generalization of Theorem~\ref{main} is possible also remains an open problem.

\subsection{Asymptotic spectrum of graphs}\label{subsect:as}

\medskip
\par\noindent
The categorical (also called direct or tensor) product $F\times G$ of two graphs $F$ and $G$ is defined as follows.
$$V(F\times G)=V(F)\times V(G)$$
and
$$E(F\times G)=\{\{(f,g),(f',g'): \{f,f'\}\in E(F)\ {\rm and}\ \{g,g'\}\in E(G)\}.$$

\medskip
\par\noindent
Let $F\to G$ denote that there exists a homomorphism of graph $F$ to graph $G$, which means an edge-preserving map from $V(F)$ to $V(G)$. It is easy to see that if $p(G)$ is a graph parameter (a function from graphs to the real numbers) which is homomorphism monotone increasing, that is, for which $F\to G$ implies $p(F)\le p(G)$, then $$p(F\times G)\le\min\{p(F),p(G)\}.$$ We say that the parameter $p$ satisfies the Hedetniemi-type equality, if the previous inequality holds with equality for every pair of graphs $F,G$. The name refers to Hedetniemi's conjecture which stated that the chromatic number is such a parameter. Although this was refuted after more than half a century (see \cite{Shitov} for the first and \cite{Tard4} for the strongest counterexample to date), some interesting graph parameters do satisfy this equality nontrivially. These include the fractional chromatic number and the complementary Lov\'asz theta number by the results in \cite{Zhufrac} and \cite{GRSS}, respectively. Let ${\cal C}_{\rm Hedet}$ denote the family of all graph parameters that satisfy the Hedetniemi-type equality.

\medskip
\par\noindent
Note that the (generalized) Mycielski construction is closely related to the categorical product in the sense that $M_r(G)$ can be defined (as it is done in fact in \cite{CTcone}, for example) by taking the categorical product of $G$ with an $r$-length path on vertices $0,1,\dots,r$ that has a loop at vertex $0$ and then identify all vertices of the form $(v,r)$ (where $v\in V(G)$).
Let ${\cal C}_{\rm Myc}$ denote the family of graph parameters $p$ for which there exists some function $g: \mathbb{R}\to \mathbb{R}$ satisfying $p(M(G))=g(p(G))$.
We can also define similarly ${\cal C}_{\rm genMyc}$ as the family of graph parameters
for which there exists a function $g_r$ for every $r$ such that $p(M_r(G))=g_r(p(G))$. Obviously ${\cal C}_{\rm genMyc}\subseteq {\cal C}_{\rm Myc}$ directly follows from the definitions. The already quoted result of Larsen-Propp-Ullman \cite{LPU} gives that the fractional chromatic number belongs to ${\cal C}_{\rm Myc}$ while our Theorem~\ref{main} gives that $\thet(G)$ also belongs to it. Tardif's result in \cite{CTcone} also shows that $\chi_f(G)$ is also a member of ${\cal C}_{\rm genMyc}$, while we do not know whether the analogous statement is true for $\thet(G)$.

\medskip
\par\noindent
Using a theory of Strassen developed to investigate the complexity of matrix multiplication Zuiddam \cite{Zuiddam} defines (in a complementary language) what he calls the asymptotic spectrum of graphs that we will denote by ${\cal AS}$. A graph parameter $p$ belongs to this family of functions if it satisfies the following four requirements.
\begin{enumerate}

\item
  Homomorphism monotonicity: $G\to H\Rightarrow p(G)\le p(H)$;

\item
  multiplicativity for the OR-product: $\forall G,H:\ p(G\cdot H)=p(G)p(H)$;

\item
  additivity for complete join: $\forall G,H:\ p(G\oplus H)=p(G)+p(H),$ where
  $V(G\oplus H)$ is the disjoint union of $V(G)$ and $V(H)$ and $E(G\oplus H)=E(G)\cup E(H)\cup\{\{g,h\}: g\in V(G), h\in V(H)\}$;

\item
  normality: $p(K_1)=1$.
\end {enumerate}

\medskip
\par\noindent
It is straightforward to prove from the above properties that if $p\in {\cal AS}$ then $C_{\rm OR}(G)\le p(G)$ holds for all graphs $G$. Using Strassen's theory (for more details on which cf. \cite{WZ}) Zuiddam \cite{Zuiddam} proved the surprising fact that $$C_{\rm OR}(G)=\min_{p\in {\cal AS}} p(G)$$ holds for all graphs $G$. This would be trivial if $C_{\rm OR}(G)$ itself would belong to ${\cal AS}$ but by well-known results of Haemers \cite{Haemers} and Alon \cite{Alon} $C_{\rm OR}(G)$ is neither multiplicative under the OR-product nor additive under the complete join. On the other hand, both $\chi_f(G)$ and $\thet(G)$ do belong to ${\cal AS}$.

\medskip
\par\noindent
Thus ${\cal C}_{\rm Hedet}, {\cal C}_{\rm Myc}$ and ${\cal AS}$ are three classes of graph parameters whose intersection contains both the fractional chromatic number and the complementary Lov\'asz theta number. Our main goal in describing these classes together is to ask whether there could be any stronger relationship among them. (Zuiddam \cite{Zuiddam} lists all known members of ${\cal  AS}$ that are of not too many types.) Note that no two of these classes are equal. It is observed in \cite{shanhed} that Zuiddam's result implies that if we had ${\cal AS}\subseteq {\cal C}_{\rm Hedet}$ then it would imply $C_{\rm OR}(G)\in {\cal C}_{\rm Hedet}$, too, therefore we must have $${\cal AS}\neq {\cal C}_{\rm Hedet}.$$ A more direct example of this non-equality is provided though by the clique number that belongs to both ${\cal C}_{\rm Hedet}$ and ${\cal C}_{\rm Myc}$ but not to ${\cal AS}$.
Since the chromatic number is neither in ${\cal AS}$ nor in ${\cal C}_{\rm Hedet}$ (an easy fact for ${\cal AS}$, since already $\chi(C_5^2)<9=[\chi(C_5)]^2$, while a highly nontrivial result refuting a half-century old conjecture in case of ${\cal C}_{\rm Hedet}$), but it does belong to ${\cal C}_{\rm Myc}$, we know that ${\cal C}_{\rm Myc}$ is also not contained in ${\cal C}_{\rm Hedet}$.  (The situation is different with ${\cal C}_{\rm genMyc}$, however, which does not contain the chromatic number either. At the same time it does trivially contain the clique number.) It might be a farfetched idea but since the parameters in ${\cal AS}$ are somehow ``the nice ones'', it would sound us to be reasonable to believe that perhaps ${\cal AS}\subseteq {\cal C}_{\rm Myc}$, for example. This, however, remains the subject of further research.

\section{Open problems}

Here we collect the open questions that (at least implicitly) came up in the above discussions and we find the most interesting.

\begin{enumerate}

\item
  Is it possible that $C_{\rm OR}(M(K_3))=4\cos{\left(\frac{2\pi}{9}\right)},$ that is that $\thet$ gives a tight upper bound also in this case as for $M(K_2)$?

\item
  Is the value of $\thet(M_r(G))$ determined by $\thet(G)$ also for $r>2$?

\item
  Do we have $C_{\rm Sp}(M_r(T_2))>2$ also for $r\ge 4$?

\item
  Does $C_{\rm OR}(M_r(G))>C_{\rm OR}(G)$ hold for every $G$ and $r$? In particular, does it hold for $r=2$ and every graph $G$?

\item
Is $${\cal AS}\subseteq {\cal C}_{\rm Myc}$$ true?
Or any other nontrivial relation among the classes of graph parameters mentioned in Subsection~\ref{subsect:as}?

  \bigskip
  \par\noindent
  Finally we also mention the following problem that is known to be equivalent to a famous open problem of Erd\H{o}s about Ramsey numbers. To formulate the problem we use the notation $M^t(G)$ to denote the $t$ times iterated Mycielskian of $G$.

\item
  Is the limit $\lim_{t\to\infty}C_{\rm OR}(M^t(K_2))$ finite?

  \smallskip
  \par\noindent
  The existence of the limit follows immediately from the fact that $M^{t-1}(G)$ is a subgraph of $M^t(G)$, therefore $C_{\rm OR}(M^{t-1}(G))\le C_{\rm OR}(M^t(G))$ for any graph $G$. For the equivalence with Erd\H{o}s's problem cf. \cite{coletrips} and the references therein.

\section{Acknowledgements}
We thank Anna Gujgiczer for useful discussions.

\end{enumerate}

\newpage

\par\noindent
{\Large{\bf A1. Appendix 1}}

\bigskip
\par\noindent
In this Appendix we give the details of the calculation leading to the cubic equation~\ref{eq:cubic} given in Section~\ref{sect:theta}

\medskip
\par\noindent
We want to calculate the solutions of the system of nonlinear equations below.
\begin{equation}
\centering
\left\{\begin{split}
\alpha^2 + x^2 &= 1 \\
\beta^2 + y^2 &= 1 \\
 \frac{1}{t-1}\alpha^2 - x^2 &=\frac{1}{\hat{t}-1}\\
\frac{1}{t-1}\alpha\beta + xy &=\frac{1}{\hat{t}-1}\\
y  &=\frac{1}{\hat{t}-1}.
\end{split}\right.
 \end{equation}
 In order to simplify the calculation, we use the $v:=t-1$ and $w:=\frac{1}{\hat{t}-1}$ substitutions. So we have
 \begin{equation}
\centering
\left\{\begin{split}
\alpha^2 + x^2 &= 1 \\
\beta^2 + w^2 &= 1 \\
 \frac{1}{v}\alpha^2 - x^2 &=w\\
\frac{1}{v}\alpha\beta + xw &=w.\\
\end{split}\right.
\end{equation}
From the first and the second equations we have $\alpha = \sqrt{1-x^2}$ and $\beta=\sqrt{1-w^2}$. Hence the third equation becomes
\[
\frac{1}{v}(1-x^2) -x^2 =w \Rightarrow x = \sqrt{\frac{1-vw}{v+1}}.
\]
Now we need to substitute the known parameters to the fourth equation instead of $x$ in order to simplify the calculation.
\begin{gather*}
\frac{1}{v}\sqrt{1-x^2}\sqrt{1-w^2} + wx = w \qquad /\cdot v, -wx\\
\sqrt{1-x^2}\sqrt{1-w^2} = vw\left(1-x\right)  \qquad /()^2\\
\left(1-x^2\right)(1-w^2)=v^2w^2\left(1-x\right)^2\\
\left(1-x^2\right)(1-w^2)=v^2w^2\left(1-2x+x^2\right) \qquad /-v^2w^2,-v^2w^2x^2\\
\left(1-x^2\right)(1-w^2) -v^2w^2 - v^2w^2x^2 = -2v^2w^2x \qquad /()^2\\
\left(1-x^2\right)^2(1-w^2)^2 + v^4w^4 + v^4w^4\left(x^2\right)^2 -2v^2w^2\left(1-x^2\right)(1-w^2) \\
- 2\left(1-x^2\right)(1-w^2)v^2w^2x^2 + 2v^4w^4x^2 = 4v^4w^4x^2
\end{gather*}
After the products and squares are calculated, $x$ is substituted with $ \sqrt{\frac{1-vw}{v+1}}$.
\begin{gather*}
\left(\frac{1-vw}{v+1}\right)^2v^4w^4 - 2\left(\frac{1-vw}{v+1}\right)^2v^2w^4 + 2\left(\frac{1-vw}{v+1}\right)^2v^2w^2 + \left(\frac{1-vw}{v+1}\right)^2w^4 \\
- 2\left(\frac{1-vw}{v+1}\right)^2w^2 + \left(\frac{1-vw}{v+1}\right)^2 - 2\frac{1-vw}{v+1}v^4w^4 - 2\frac{1-vw}{v+1}w^4 + 4\frac{1-vw}{v+1}w^2 \\
- 2\frac{1-vw}{v+1} + v^4w^4 + 2v^2w^4 - 2v^2w^2 + w^4 - 2w^2 + 1 = 0; \qquad /\cdot(v+1)^2\\
(1-vw)^2v^4w^4 - 2(1-vw)^2v^2w^4 + 2(1-vw)^2v^2w^2 + (1-vw)^2w^4- 2(1-vw)^2w^2 \\
+ (1-vw)^2 - 2(1+v)(1-vw)v^4w^4 - 2(1+v)(1-vw)w^4 + 4(1+v)(1-vw)w^2 \\
- 2(1+v)(1-vw) + (1+v)^2v^4w^4 + 2(1+v)^2v^2w^4 - 2(1+v)^2v^2w^2 + (1+v)^2w^4 \\
- 2(1+v)^2w^2 + (1+v)^2 = 0;\\
v^6w^6 + 2v^6w^5 + v^6w^4 - 2v^4w^6 + 4v^4w^4 - 2v^4w^2 + 4v^3w^5 + 4v^3w^4 - 4v^3w^3 \\
- 4v^3w^2 + v^2w^6 + 2v^2w^5 - v^2w^4 - 4v^2w^3 - v^2w^2 + 2v^2w + v^2 = 0; \qquad /:v^2\\
v^4w^6 + 2v^4w^5 + v^4w^4 - 2v^2w^6 + 4v^2w^4 - 2v^2w^2 + 4vw^5 + 4vw^4 - 4vw^3 - 4vw^2 \\
+ w^6 + 2w^5 - w^4 - 4w^3 - w^2 + 2w + 1 = 0;
\end{gather*}
Fortunately, the last polynomial can be factorized as follows
\begin{gather*}
    (w + 1)^2 (v w + w - 1) (v^3 w^3 - v^2 w^3 + v^2 w^2 - v w^3 + 2 v w^2 - v w + w^3 + w^2 - w - 1) = 0.
\end{gather*}
Hence we have three factor polynomials, thus the roots of the original polynomial are all the roots of the following polynomials, $p_1(w):=w+1,p_2(w,v)=vw+w-1,p_3(v,w):=v^3 w^3 - v^2 w^3 + v^2 w^2 - v w^3 + 2 v w^2 - v w + w^3 + w^2 - w - 1$, From the definition of the first polynomial $w=-1$, $\hat{t} = 0$ so it is not the solution which is expected. And $p_2(v,w) = 0$ shows that $w = \frac{1}{v+1}$, therefore $\hat{t}=t+1$. We will see that $t+1$ is never smaller than the maximal root of the third polynomial. (It is in fact not hard to see that the root $\hat{t}=t+1$ belongs to the ``degenerate'' strict vector coloring of $M(G)$ in which we have $u_i^*=v_i^*$ for all $i$. It can also be interpreted as an, in that case optimal, strict vector coloring of the ``generalized'' Mycielskian $M_1(G)$ of $G$, that is, when only one vertex $z_{M_1(G)}$ is added to $G$ that is made adjacent to all original vertices.)
Thirdly we have the equation
\begin{gather*}
    v^3 w^3 - v^2 w^3 + v^2 w^2 - v w^3 + 2 v w^2 - v w + w^3 + w^2 - w - 1 = 0;\\
    \frac{(t-1)^3}{(\hat{t}-1)^3} - \frac{(t-1)^2}{(\hat{t}-1)^3} + \frac{(t-1)^2}{(\hat{t}-1)^2} -  \frac{t-1}{(\hat{t}-1)^3} + 2  \frac{t-1}{(\hat{t}-1)^2} -  \frac{t-1}{\hat{t}-1} + \frac{1}{(\hat{t}-1)^3} \\
    + \frac{1}{(\hat{t}-1)^2} - \frac{1}{\hat{t}-1} - 1 = 0; \qquad /\cdot (\hat{t}-1)^3\\
 (t-1)^3  - (t-1)^2  + (t-1)^2 (\hat{t}-1) - (t-1)  + 2 (t-1) (\hat{t}-1) - (t-1) (\hat{t}-1)^2 \\
 + 1 + (\hat{t}-1) - (\hat{t}-1)^2 - (\hat{t}-1)^3 = 0;\\
  t^3 + t^2\hat{t} - 5t^2 - t\hat{t}^2 + 2t\hat{t} + 3t - \hat{t}^3 + 3\hat{t}^2 - 3\hat{t} + 1 = 0;\\
  \hat{t}^3 + (t-3)\hat{t}^2 + (3-2t-t^2)\hat{t} + (-t^3 +5t^2 -3t-1) = 0.
\end{gather*}

\medskip
\par\noindent
The last equation is already the cubic equation we had as equation~(\ref{eq:cubic}) in Sextion~\ref{sect:theta}.

\medskip
\par\noindent
It remains to check that $\frac43t\cos{\left(\frac13{\rm arccos}\left(1 - \frac{27}{4t}+\frac{27}{4t^2}\right)\right)}-\frac13t+1 \leq t+1$ that we stated on the way. After rearrangement we have $\cos{\left(\frac13{\rm arccos}\left(1 - \frac{27}{4t}+\frac{27}{4t^2}\right)\right)}\leq 1$ and it is a well-known property.\\
\vspace{0.3 cm}

\bigskip
\bigskip
\bigskip

\par\noindent
{\Large{\bf A2. Appendix 2}}

\bigskip
\par\noindent
In this Appendix we prove that $x_1, x_2\le 1$, where $x_1$ and $x_2$ are the values obtained from Equation~(\ref{eq:star}) when substituting $k=1,2$, respectively.

\medskip
\par\noindent
So we have to see that we indeed have  \[
\frac43t\cos{\left(\frac13{\rm arccos}\left(1 - \frac{27}{4t}+\frac{27}{4t^2}\right)-\frac{2\pi k}{3}\right)}-\frac13t+1 \leq 1\hspace{0.5 cm} \text{for }k=1,2.
\]

Subtracting $1$ from both sides and then dividing by $\frac43t$ we arrive to the equivalent inequality
$$\cos{\left(\frac13{\rm arccos}\left(1 - \frac{27}{4t}+\frac{27}{4t^2}\right)-\frac{2\pi k}{3}\right)}-\frac14 \leq 0.$$

First consider the $k=1$ case. Writing simply $X$ in place of $1 - \frac{27}{4t}+\frac{27}{4t^2}$
some equivalent manipulations give us

$$\sin\left(\frac{\pi}{2}-\frac13{\rm arccos}(X)+\frac{2\pi}{3}\right)-\frac14\le 0,$$

which is further equivalent to

$$\sin\left(\frac{\pi}{2}-\frac13\left(\frac{\pi}{2}-{\rm acsin}(X)\right)+\frac{2\pi}{3}\right)-\frac14\le 0$$

that is the same as

$$\sin\left(\frac13{\rm arcsin}(X)+\pi\right)-\frac14\le 0$$

that is

$$-\sin\left(\frac13{\rm arcsin}(X)\right)-\frac14\le 0,$$

so (putting back $X$) we need to verify that

$$\sin{\left(\frac13{\rm arcsin}\left(1-\frac{27}{4t}+\frac{27}{4t^2}\right)\right)} \geq -\frac14.$$
The minimum value of $1-\frac{27}{4t}+\frac{27}{4t^2}$ is $-\frac{11}{16}$, so the inequality holds for all $1\leq t$.

In case of $k=2$ we get to the following calculation by equivalent manipulations.

\begin{gather*}
    \frac43t\cos{\left(\frac13{\rm arccos}\left(1-\frac{27}{4t}+\frac{27}{4t^2}\right)-\frac{4\pi}{3}\right)} - \frac13t + 1 \leq - \frac13t + 1\\
   \cos{\left(\frac13{\rm arccos}\left(1-\frac{27}{4t}+\frac{27}{4t^2}\right)-\frac{4\pi}{3}\right)} \leq 0\\
   \sin{\left(\frac{\pi}{2}-\frac13{\rm arccos}\left(1-\frac{27}{4t}+\frac{27}{4t^2}\right)+\frac{4\pi}{3}\right)} \leq 0\\
   \sin{\left(\frac{11\pi}{6}-\frac13{\rm arccos}\left(1-\frac{27}{4t}+\frac{27}{4t^2}\right)\right)} \leq 0\\
   \sin{\left(-\frac{\pi}{6}-\frac13{\rm arccos}\left(1-\frac{27}{4t}+\frac{27}{4t^2}\right)\right)} \leq 0\\
\end{gather*}

For the above it is enough to have

\begin{align*}
    0 \leq \frac{\pi}{6}+\frac13{\rm arccos}\left(1-\frac{27}{4t}+\frac{27}{4t^2}\right)
\end{align*}

and since the second term here is at least $-\frac{\pi}{2}$, this inequality holds for any $t>0$.

\bigskip
\bigskip
\bigskip

\par\noindent
{\Large{\bf A3. Appendix 3}}

\bigskip
\par\noindent
In this Appendix we prove Lemma~\ref{lem:maradt} that stated the following inequalities.

    \begin{enumerate}
    \item $0< \eta$;
    \item $\frac{(m(t)-1)^2}{t-1} \leq \gamma^2$;
    \item $0< \delta$.
    \end{enumerate}

\proof
We have $\left(\frac{m(t)-1}{t-1}\right)^2 \leq \frac{(m(t)-1)^2}{t-1}$ by $t\ge 2$. Therefore $\frac{(m(t)-1)^2}{t-1} \leq \gamma^2$ implies $\gamma^2\geq \left(\frac{m(t)-1}{t-1}\right)^2$ which is equivalent to $$0<\gamma\frac{t-1}{m(t)-1}-\frac{1}{\gamma\frac{t-1}{\hat{t-1}}}=\delta$$
by $m(t)>t$ and $\gamma>0$ and $t\ge 2$.

So $0< \delta$ follows if we know $\frac{(m(t)-1)^2}{t-1} \leq \gamma^2$ and even strict inequality follows from $t>2$.


\medskip
\par\noindent
By the formula we obtained for $\eta$ we need to show
\[
\frac{1}{\delta(m(t)-1)}\left(\gamma^2\delta + \gamma^3\frac{1}{m(t)-1}-\gamma^3\right)>0
\]

\par\noindent
We show first that $$\gamma^2\delta + \gamma^3\frac{1}{m(t)-1}-\gamma^3>0,$$ then the whole inequality will follow if we also show $\delta>0$, since $m(t)>1$ is true. Substituting the formula $\delta = \gamma\frac{t-1}{m(t)-1} - \frac{1}{\gamma\frac{t-1}{m(t)-1}}$ we have for $\delta$, we obtain after rearrangement and multiplication with $m(t)-1>0$ that the claimed inequality is equivalent to
\begin{equation}\label{eq:gammasq}
\gamma^2(t-m(t)+1)\ge \frac{(m(t)-1)^2}{t-1},
\end{equation}

\smallskip
\par\noindent
Substituting our formula obtained for $\gamma^2$ this turns out to be equivalent to $-\frac{(m(t) - 1)^2 (t - m(t) + 1)^3}{2 (t - 1) t (m(t) - 2) (t-m(t))} \geq \frac{(m(t)-1)^2}{t-1}$.
After further rearrangement this is seen to be equivalent to
\[
\frac{\left(\frac{(1+t-m(t))^3}{t(m(t)-t)(m(t)-2)}-2\right)(m(t)-1)^2}{2(t-1)} \geq 0, \]
and after multiplying by $2(t-1)>0$ and dividing by $(m(t)-1)^2>0$ we are left with the equivalent formula
$$\frac{(1+t-m(t))^3}{t(m(t)-t)(m(t)-2)}-2 \geq 0.$$
Some further equivalent manipulations follow to obtain an inequality we can directly see to be true:
We multiply both sides by $t(m(t)-t)(m(t)-2)>0$ and simplify the inequality to obtain
$$-m(t)^3 + (t+3)m(t)^2 + (-t^2-2t-3)m(t) + (t^3 - t^2 +3t+1)\geq 0.$$
We know that $t,m(t)$ satisfy the equation\\
$m(t)^3 + (t-3)m(t)^2 + (3-2t-t^2)m(t) + (-t^3 +5t^2 -3t-1) = 0$, therefore we can write
\begin{gather*}
-m(t)^3 + (t+3)m(t)^2 + (-t^2-2t-3)m(t) + (t^3 - t^2 +3t+1) =\\
-[m(t)^3 + (t-3)m(t)^2 + (3-2t-t^2)m(t) + (-t^3 +5t^2 -3t-1)] - 2 t^2 m(t) + 4 t^2 +\\
+2 t m(t)^2 - 4 t m(t) = -2 t (m(t) - 2) (t - m(t)),
\end{gather*}
since $t < m(t)$, the previous inequality follows as $-2 t (m(t) - 2) (t - m(t))\geq 0$ is true.\\
Thus we, in particular, proved inequality (\ref{eq:gammasq}) above and since $t-m(t)+1 \leq 1$, this also proves
\begin{equation}\label{eq:gam}
\frac{(m(t)-1)^2}{t-1} \leq \gamma^2.
\end{equation}
\medskip
\par\noindent
What is left is to see that $\delta>0$ also holds. Using $\delta = \gamma\frac{t-1}{m(t)-1}-\frac{1}{\gamma\frac{t-1}{\hat{t-1}}}$ and multiplying by the denominator of the second term (that we know to be positive) here, the required inequality turns out to be equivalent to
$$\gamma^2\geq \left(\frac{m(t)-1}{t-1}\right)^2.$$
The truth of this follows from $\left(\frac{m(t)-1}{t-1}\right)^2 \leq \frac{(m(t)-1)^2}{t-1}$ (implied by $t\ge 2$) and inequality (\ref{eq:gam}) above.
This completes the proof of Lemma~\ref{lem:maradt}
\qed

\bigskip
\bigskip
\bigskip

\par\noindent
{\Large{\bf A4. Appendix 4}}

\bigskip
\par\noindent
In this appendix, we will calculate the simpler form stated in Section~\ref{sect:theta} of the following equation.
\begin{gather*}
\gamma\delta + \gamma^2\frac{1}{m(t)-1}-\gamma^2-\gamma^2\frac{1}{m(t)-1}-\frac{1}{m(t)-1}\gamma\delta-\gamma^2\frac{1}{(m(t)-1)^2} + \gamma^2\frac{1}{m(t)-1}=\\
-\ \frac{1}{\delta(m(t)-1)}\left(\gamma^2\delta + \gamma^3\frac{1}{m(t)-1}-\gamma^3\right) -1
\end{gather*}
Similarly to the calculation in Appendix 1 we use the $v:=t-1,w:=\frac{1}{m(t)-1}$ substitutions. According to these notations we have the following expressions for our parameters:
\begin{align*}
    \delta =& \gamma vw - \frac{1}{\gamma vw},\\
    \eta =& \frac{w}{\delta}(\gamma^2\delta + \gamma^2w-\gamma^3),
\end{align*}

hence we need to solve the following equation:
\begin{gather*}
    \gamma\delta + \gamma^2w -\gamma^2 -\gamma^2w -\gamma\delta w -\gamma^2w^2 + \gamma^2w = -\frac{w}{\delta}(\gamma^2\delta +\gamma^3w-\gamma^3)-1;\qquad /\cdot \delta\\
    \gamma\delta^2 -\gamma^2\delta -\gamma w\delta^2  -\gamma^2w^2\delta + \gamma^2w\delta = -\gamma^2w\delta -\gamma^3w^2+\gamma^3w-\delta;\\
    \gamma\delta^2 -\gamma^2\delta -\gamma w\delta^2  -\gamma^2w^2\delta + \gamma^2w\delta  +\gamma^2 w\delta +\gamma^3w^2-\gamma^3w+\delta=0;\\
    \gamma\left(\gamma vw - \frac{1}{\gamma vw}\right)^2 -\gamma^2\left(\gamma vw - \frac{1}{\gamma vw}\right) -\gamma w\left(\gamma vw - \frac{1}{\gamma vw}\right)^2  -\gamma^2w^2\left(\gamma vw - \frac{1}{\gamma vw}\right) \\
    + \gamma^2w\left(\gamma vw - \frac{1}{\gamma vw}\right)
    +\gamma^2 w\left(\gamma vw - \frac{1}{\gamma vw}\right) +\gamma^3w^2-\gamma^3w+\left(\gamma vw - \frac{1}{\gamma vw}\right)=0.
\end{gather*}
Multiplying both sides by $\gamma^2 v^2w^2$ we obtain
\begin{gather*}
        \gamma(\gamma^2 v^2w^2 - 1)^2 -\gamma^2(\gamma^3 v^3w^3 - \gamma vw) -\gamma w(\gamma^2 v^2w^2 - 1)^2  -\gamma^2w^2(\gamma^3 v^3w^3 - \gamma vw) \\
        + \gamma^2w(\gamma^3 v^3w^3 - \gamma vw)  +\gamma^2 w(\gamma^3 v^3w^3 - \gamma vw) +\gamma^5v^2w^4-\gamma^5v^2w^3+(\gamma^3 v^3w^3 - \gamma vw)=0.\\
        -\gamma^5v^4w^5 + \gamma^5v^4w^4 - \gamma^5v^3w^5 + 2\gamma^5v^3w^4 - \gamma^5v^3w^3 + \gamma^5v^2w^4 - \gamma^5v^2w^3 + \gamma^3v^3w^3 \\
        + 2\gamma^3v^2w^3 - 2\gamma^3v^2w^2 + \gamma^3vw^3 - 2\gamma^3vw^2 + \gamma^3vw - \gamma vw - \gamma w + \gamma=0;\qquad /:\gamma\\
        (-v^4w^5 + v^4w^4 - v^3w^5 + 2v^3w^4 - v^3w^3 + v^2w^4 - v^2w^3)\gamma^4 \\
        + (v^3w^3 + 2v^2w^3 - 2v^2w^2 + vw^3 - 2vw^2 + vw)\gamma^2 - vw - w + 1 = 0.\\
\end{gather*}
So we have the (in $\gamma^2$) quadratic equation $a\gamma^4 + b\gamma^2 + c = 0$ from where the proof continued in Section~\ref{sect:theta} of the paper.

\bigskip
\bigskip
\bigskip

\end{document}